\documentstyle[11pt, amssymb,amsfonts,amscd, amsmath, xypic]{article}

\xyoption{all}
\CompileMatrices

\newcommand{\triv}{{\sf{triv}}}
\newcommand{\ZZ}{{\sf{Z}}}
\newcommand{\nc}{\newcommand}
\newcommand{\qed}{$\enspace\square$}
\newcommand{\Vx}{{\mathcal{V}}}

\newcommand{\dis}{\displaystyle}

\newcommand{\anol}{{A_0}}
\newcommand{\fin}{^{^{\sf{fin}}}}

\newtheorem{theorem}[equation]{Theorem}
\newtheorem{proposition}[equation]{Proposition}
\newtheorem{corollary}[equation]{Corollary}
\newtheorem{definition}[equation]{Definition}

\newtheorem{lemma}[equation]{Lemma}
\newtheorem{example}[equation]{Example}

\newenvironment{case}{\left\{\begin{array}{rl}}{\end{array}\right.}

\nc{\dlim}{\mathop{\lim\limits_{\longrightarrow}}}
\nc{\FB}{{\mathcal F\mathcal B}}
\nc{\B}{{\mathcal B}}
\nc{\BS}{{\underline{\B}}}
\nc{\Brat}{{\mathcal B_{\mathrm frac}}}
\nc{\Grad}[1]{\mathsf{Gr}^{\mathsf{ad}}_{#1}}
\nc{\vdim}{{\mathrm v.}\dim}
\nc{\Qu}{{Q}}
\nc{\FM}{{\mathfrak M}}
\nc{\TM}{{\widetilde M}}
\nc{\BP}{{\mathsf{P}}}
\nc{\QQ}{{\mathsf{Q}}}
\nc{\PPP}{{\PP^1\times_\tau\PP^1}}
\nc{\iz}{{i_z}}
\nc{\iw}{{i_w}}
\nc{\Quot}{{{\mathrm Quot}_\CO^0}}
\nc{\BQuot}{{{\mathrm Quot}_\B}}
\nc{\Subs}{{{\mathrm Sub}_\CO^0}}
\nc{\BSubs}{{{\mathrm Sub}_\B}}
\nc{\CSubs}{{{\mathrm Sub}_\C}}
\nc{\CM}{{\mathcal M}}
\nc{\epss}{{\eps^{-1}}}
\nc{\CG}{{\C\Gamma}}
\nc{\Diff}{{\mathop{\mathrm Diff}}}
\nc{\PGL}{{\mbox{PGL}}}
\nc{\longleftrightarrows}{\quad\rlap{\raisebox{-2pt}{$\longleftarrow$}}%
\raisebox{2pt}{$\longrightarrow$}\quad}
\nc{\G}{\Gamma}
\nc{\bplus}{\mbox{$\bigoplus$}}
\nc{\btimes}{\mbox{$\bigotimes$}}

\def\beq{\begin{equation}}
\def\proj{{\sf{Proj}}}
\def\eeq{\end{equation}}

\newcommand{\iso}{{\;\stackrel{\sim}{\longrightarrow}\;}}
\newcommand{\cd}{\!\cdot\!}
\newcommand{\vo}{${\sf {(0)}}\;$}
\newcommand{\vi}{${\sf {(i)}}\;$}
\newcommand{\vii}{${\sf {(ii)}}\;$}
\newcommand{\viii}{${\sf {(iii)}}\;$}

\newcommand{\id}{{\mbox{\sl Id}}}

\newcommand{\into}{\,\hookrightarrow\,}
\newcommand{\too}{\,\longrightarrow\,}
\newcommand{\onto}{\,\twoheadrightarrow\,}

\newcommand{\Spec}{{\mbox{\sf Spec}^{\,}}}

\newcommand{\Hom}{{\mbox{\sl Hom}}}

\newcommand{\lef}{_{_{\sf{left}}}}
\newcommand{\om}{\omega}
\newcommand{\op}{\mbox}

\newcommand{\eps}{{\epsilon}}

\newcommand{\h}{{\mathfrak{h}}}

\newcommand{\X}{{\mathcal{X}}}
\def\C{{\mathbb{C}}}

\def\PP{{\mathbb{P}}}
\def\bp{{\mathbf{P}}}
\def\Z{{\mathbb{Z}}}

\def\K{{\mathcal K}}

\def\CO{{\mathcal O}}

\def\A{{\mathcal{A}}}

\def\gln{{\mathfrak{g}\mathfrak{l}}_n(\C)}

\def\bbm{{\mathbb{M}}}
\def\bcm{{\boldsymbol{\mathcal{M}}}}

\def\ccirc{{}_{^\circ}}

\newcommand{\modd}{\mathop{\op{\sf mod}}}
\newcommand{\coh}{\mathop{\op{\sf coh}}}
\newcommand{\qgr}{\mathop{\op{\sf qgr}}}
\newcommand{\gr}{\mathop{\op{\sf gr}}}
\newcommand{\tor}{\mathop{\op{\sf tor}}}

\newcommand{\lra}{{\longrightarrow}}

\nc{\Fl}{\mathop{\rm Fl}\nolimits}

\nc{\newi}[1]{\noindent{\rm(#1)}}
\nc{\CC}{{\mathcal{C}}}
\nc{\CE}{{\mathcal{E}}}
\nc{\CF}{{\mathcal{F}}}
\nc{\CH}{{\mathcal{H}}}
\nc{\Tor}{\mathop{{\rm Tor}}\nolimits}
\nc{\CHom}{{\underline{{\mathcal{H}}om}}}
\nc{\CExt}{{\underline{{\mathcal{E}}xt}}}
\renewcommand{\Hom}{\mathop{{\rm Hom}}\nolimits}
\nc{\Ext}{\mathop{{\rm Ext}}\nolimits}
\nc{\End}{\mathop{{\rm End}}\nolimits}
\nc{\Ker}{\mathop{\mathbf{Ker}}\nolimits}
\renewcommand{\Im}{\mathop{\mathbf{Im}}\nolimits}
\nc{\Coker}{\mathop{\mathbf{Coker}}\nolimits}
\nc{\k}{{\mathsf{k}}}
\nc{\TV}{{\widetilde V}}
\nc{\CN}{{\mathcal N}}

\nc{\gtimes}{\otimes_{\!_\CG}}
\nc{\Rep}{\mathop{\mathbf{Rep}}\nolimits}
\nc{\CU}{{\mathcal U}}

\nc{\CR}{{\mathcal R}}
\nc{\ip}[1]{\left[{#1}\right]}
\nc{\mh}{{\tilde h}}
\nc{\ppg}{{\PP^2_{\small \Gamma}}}
\nc{\ppgt}[1]{{\PP^2_{{\small \Gamma}}}}
\nc{\plg}{{\PP^1_{{\small \Gamma}}}}

\begin{document}
\setlength{\parindent}{6mm}
\setlength{\parskip}{3pt plus 5pt minus 0pt}

\centerline{\Large {\bf Quiver varieties and a noncommutative $\PP^2$}}

\vskip 6mm
\centerline{\large {\sc {V. Baranovsky, V. Ginzburg and A. Kuznetsov}}}
\vskip 6pt

\begin{abstract}
{\footnotesize 
To any  finite group $\Gamma \subset SL_2(\C)$ and
each  element $\tau$ in the center of the
group algebra of $\G$, we associate a category,
${\mathcal{C}}\!oh(\PP^2_{_{\!\Gamma^{\!},^{\!}\tau}},\PP^1).$ 
It is defined as a suitable quotient of the category
of graded modules over 
(a graded version of) the {\it deformed
preprojective algebra}
introduced by Crawley-Boevey and Holland.
The category
${\mathcal{C}}\!oh(\PP^2_{_{\!\Gamma^{\!},^{\!}\tau}},\PP^1)$
  should be
thought of as the category of coherent sheaves on a
`noncommutative projective space', $\PP^2_{_{\!\Gamma^{\!},^{\!}\tau}},$ equipped
with a framing at  $\PP^1$, the line at infinity.
We
 establish an isomorphism between the moduli space of
\textit{torsion free} objects
of ${\mathcal{C}}\!oh(\PP^2_{_{\!\Gamma^{\!},^{\!}\tau}},\PP^1)$ and
the Nakajima quiver
variety arising from $\G$ via the  McKay correspodence.

We  apply the above isomorphism
 to deduce a generalization of the Crawley-Boevey and Holland
conjecture,
saying that the moduli space of `rank 1' projective
modules over the deformed
preprojective algebra is isomorphic to a particular
quiver variety. This reduces, for $\G=\{1\}$, to the recently
obtained parametrisation of the isomorphism classes of right ideals
in the first Weyl algebra, $\sf{A_1}$, by points of the
Calogero-Moser space, due to  Cannings-Holland and
Berest-Wilson. Our approach is
 algebraic and is based on 
a monadic description  of torsion free sheaves on 
$\PP^2_{_{\!\Gamma^{\!},^{\!}\tau}}$. It is
totally different from the
one used by Berest-Wilson, involving $\tau$-functions.}
\end{abstract}
\medskip

\centerline{\bf Table of Contents}
\smallskip

$\hspace{20mm}$ {\footnotesize \parbox[t]{115mm}{
1.{ $\;\,\,$} {\tt Introduction} \newline
2.{ $\;\,\,$} {\tt Sheaves on a
non-commutative {\sf{Proj}}-scheme} \newline
3.{ $\;\,\,$} {\tt Torsion free sheaves on
$\PP^2_{_{\!\Gamma^{\!},^{\!}\tau}}$} \newline
4.{ $\;\,\,$} {\tt Interpretation of  Nakajima varieties} \newline
5.{ $\;\,\,$} {\tt Proof of generalized Crawley-Boevey 
and Holland Conjecture} \newline
6.{ $\;\,\,$} {\tt Appendix A: Graded preprojective algebra} \newline
7.{ $\;\,\,$} {\tt Appendix B: Algebraic generalities} \newline
8.{ $\;\,\,$} {\tt Appendix C: Minuscule classes}
}}

\bigskip
\pagebreak[3]

\section{Introduction}
\setcounter{equation}{0}
\subsection{Noncommutative Algebraic Geometry}
A standard construction of algebraic geometry
associates to  any graded,
finitely generated commutative $\C$-algebra $A=\bigoplus_{i\geq 0}\,A_i,$
with $A_0=\C,$ a projective scheme, $\proj A$. Then
a well-known theorem
essentially due to Serre  says that the abelian category of coherent
sheaves on the scheme  $\proj A$ is equivalent to 
 $\qgr(A):=\gr(A)/\tor(A)$, the quotient of the abelian
 category  of finitely
  generated graded $A$-modules by the Serre subcategory 
of finite dimensional modules. 
Now, the category $\qgr(A)$ makes sense for any {\it non-commutative}
graded algebra $A$ as well, in which case  we let
$\gr(A)$ stand for the category of finitely generated {\it right}
$A$-modules, for definitiveness.
One may think of $\qgr(A)$
as  a category of coherent
sheaves on a `non-commutative scheme'  $\proj A$. This is
the (limited) framework in which `non-commutative  algebraic geometry'
will be understood in this paper.

 Of course, for the category $\qgr(A)$ to be reasonably
nice, the non-commutative  algebra $A$ can not be arbitrary,
and it has to satisfy some reasonable conditions. 
Throughout this paper we assume that 
$A = \bigoplus_{n \geq 0} A_n$ is a
positively graded $\C$-algebra such that all graded components,
$A_n$, are finite
dimensional over $\C$.
Write $\modd(A)$  for the category of all (not necessarily graded)
finitely
  generated right $A$-modules.

We will use the following definition, cf. [AZ]:

\begin{definition}\label{reg}\label{cond}
An algebra $A$
as above is called {\em strongly regular of dimension $d$} if

\vo $A_0$ is a (not necessarily commutative) 
semisimple $\C$-algebra.

\vi $A$ has finite global  dimension equal to $d$, that is $d$
  is the  minimal integer, such  that 
\centerline{$\Ext^{>d}_{\modd(A)}(M,N)=0,$ for
  all $M,N\in\modd(A)$.}

\vii $A$ is a Noetherian algebra of polynomial growth,
i.e., there exist integers $m,n >0$  
$\hphantom{x}\quad\;\quad$such that: 
$\dim_{_\C} A_i \leq m
\cdot i^n,$ for all $i\gg 0$;

\viii $A$  is  Gorenstein   with  parameters  $(d,l)$, i.e.
  $\Ext^i_{\modd(A)}(\anol,A)=\begin{case}\anol(l),  & \text{if  $i=d$}\\0, &
    \text{else\,.}\end{case}$
\end{definition}

\begin{remark} In [AZ],
the authors consider only the algebras $A$ such that $A_0=\C$,
and called such an $A$  {\it regular} if the
conditions \ref{reg}(i)-(iii) above hold.
It has been shown in [AZ] that any strongly
regular algebra in the sense of Definition \ref{reg} satisfies the following

{\bf{${\boldsymbol{\chi}}$-condition:}}\quad $\dis\dim_{_\C}  \Ext^i_{\modd(A)}
  (\anol, M)  < \infty,$ {\it  for  all $i\geq 0$ and }
 $M  \in \modd(A)$. \newline
To see this, note  that  any finitely
generated right $A$-module $M$  can be included in a short
exact  sequence $0 \to  M' \to  P \to  M \to  0$, where  $P$ is  a free
$A$-module  of  finite  rank,  and $M'$  is  automatically
finitely generated due to the noetherian property. Once we know the
Gorenstein property, the $(\chi)$-condition can be proved easily by
considering the long exact sequence of
$\Ext$-groups and using the descending 
induction starting at $i = {gl.  dim} (A).\;\;\lozenge$
\end{remark}\smallskip

Given a strongly regular algebra $A$,
we 
let $\pi: \gr(A) \onto \qgr(A)$ be the projection functor.
The objects of $\qgr(A)$ will be referred to as `sheaves
on $\proj A$', and we will often write $\coh(\proj A)$
instead of $\qgr(A)$.

For a graded $A$-module $M$ and $k\in\Z$, write
$M(k)$ for the same module with the grading being
shifted by $k$. For each $k\in\Z$, we
put  $\CO(k):= \pi\bigl(A(k)\bigr),$
a sheaf on $\proj A$.
Similarly, for any sheaf $E=\pi(M)$ we write $E(k)$
  for  the sheaf $\pi(M(k))$.
Further, for any sheaves $E,F\in \qgr(A)$, let $\Ext^p(E,F)$
be the $p$-th derived functor of the $\Hom$-functor:
$\Hom(E,F) = \Hom_{\qgr(A)}(E,F)$.

The crucial properties 
of the category $\coh(\proj A)=\qgr(A)$ that follow from
the strong regularity of $A$ are:

\begin{description}
\item[Ampleness, see \cite{AZ}:]  The sequence $\{\CO(i)\}_{i\in\Z}$ is
  ample, that is, for  any  $E\in\coh(\proj A)$,  there exists  an
  epimorphism: $\CO(-n)^{\oplus m}\onto E$, and for any epimorphism:
  $E\onto  F,$  the  morphism: $\Hom(\CO(-n),E)\to\Hom(\CO(-n),F)$  is
  surjective for $n\gg0$. 


\item[Serre duality, see \cite{YZ}:]  
There are integers $d\ge 0$ (dimension) and
$l$ (index of the canonical class) such that one has functorial isomorphisms
\begin{equation}\label{serre}
\Ext^i(E,F)\cong\Ext^{d-i}(F,E(-l))^\vee \quad,\quad\forall
E,F\in \coh(\proj A)
\end{equation}
where $(-)^\vee$  stands for the  dual in the category  of $\C$-vector
spaces.
\end{description}

Write $\gr\lef(A)$  for the abelian
 category  of finitely
  generated graded  {\it left} $A$-modules,
and $\pi\lef:\,
\gr\lef(A)\too  \qgr\lef(A):=\gr\lef(A)/\tor\lef(A),$ 
for the projection to the quotient  by the Serre subcategory 
of finite dimensional modules. 
Observe that the left action of the algebra $A$ on itself
by multiplication induces, for each $i$, natural
morphisms: $A_k\too \Hom_{\qgr(A)}(\CO(i), \CO(i+k))$.
This gives, for any $E\in  \qgr(A)$, a
graded {\it left} $A$-module structure on the
graded space $\oplus_{k\ge 0}\Hom(E,\CO(k))$. 
Thus, $\CHom(E,\CO)=
\pi\lef\bigl(\oplus_{k\ge 0}\Hom(E,\CO(k))\bigr)$
is a well defined object of $\qgr\lef(A)$.
This way we have defined an {\em internal Hom-functor}
$\,\CHom(-,\CO): \qgr(A)\to \qgr\lef(A)$.
Note  that it takes right modules  to left modules,
and the other way around.
The functor $\CHom(-,\CO)$ is  left exact, and we write
$\CExt^p(-,\CO): \qgr(A)\to \qgr\lef(A)$
for the corresponding derived functors.
One can check that
    \begin{equation}\label{cextp}
    \CExt^p(E,\CO) =
\pi\lef(\oplus_{k\ge0}\Ext^p(E,\CO(k)))\quad,\quad\forall p\geq 0.
    \end{equation}

For a sheaf $E\in \coh(\proj A)$, 
we define
 $\,H^p(\proj A, E)$
$:= \Ext^p(\CO,E).$ One has a 
functorial isomorphism:
 $H^0\bigl(\proj A\,,\,\CHom(E,\CO)\bigr)\cong\Hom(E,\CO),$ and also
a functorial spectral sequence:
\[ E^{p,q}_2=H^p\bigl(\proj A\,,\,\CExt^q(E,\CO)\bigr)\,
\Longrightarrow\,\Ext^\bullet(E,\CO).\]
\begin{definition}\label{free}\vi
A sheaf $E\in \qgr(A)$ is called 
{\em locally free} if $\CExt^p(E,\CO)=0\,,\,\forall p>0.$

\vii A sheaf $E\in \qgr(A)$ is called {\em torsion free} if it admits an
embedding into a locally free sheaf.
\end{definition}

The sheaves $\CO(i)$ are locally free. Moreover,
$\CHom(\CO(i),\CO(j))=\CO(j-i)$. 

\begin{remark}
It is easy to see that in the case of a smooth commutative projective
variety $X$ the definitions given above are equivalent to the standard
ones. $\quad\lozenge$
\end{remark}

\subsection{Noncommutative $\PP^2_{_\Gamma}$}
In \S3 we  will study torsion-free
sheaves on a particular 2-dimensional non-commutative
scheme analogous to $\PP^2$.
Specifically, let $(L, \om)$ be a 2-dimensional
symplectic vector space, 
and $\G\subset Sp(L)$ a finite
subgroup. 
We form the graded algebra: $TL[z]$, the polynomial algebra
    in a dummy variable $z$ (
placed in degree 1)
with coefficients in the tensor algebra of
    the vector space $L$.
 Let $(TL[z]) \# \Gamma$ be the smash product
    of $TL[z]$ with $\C\Gamma$, the group algebra of $\Gamma$, acting
    naturally on $TL$ and trivially on $z$.

To any element $\tau$ in $\ZZ(\C\G)$, the center of  $\CG$,
 we associate following
Crawley-Boevey and Holland [CBH], a
graded algebra $A^\tau$ as follows.
\begin{definition}\label{Atau}
$\quad\dis
A^\tau =  \Bigl(\bigl(TL[z]\bigr)\# \G\Bigr)
\Big/  
\big\langle
\!\big\langle u\cd v - v\cd u = \om(u,v)\cdot\tau\cdot  z^2 \big\rangle
\!\big\rangle_{u,v\in L},\;
$\newline
where $\big\langle
\!\big\langle\ldots\big\rangle
\!\big\rangle$ stands for the two-sided ideal generated by the indicated
relation.
\end{definition}

It is convenient to
 choose and fix a symplectic basis $\{x,y\}$ in $L$,
 and to identify $L$ with $\C^2$ and
$Sp(L)$ with $SL_2(\C)$. Then $TL$ gets identified
with a free associative algebra on two generators.
Writing $\C  \langle  x,y,z  \rangle$ for 
the
free associative algebra generated by $x$, $y$ and $z$,
we have

$$
A^\tau  =  \bigl(\C  \langle  x,y,z  \rangle  \#  \Gamma\bigr)  \left/  \Big\langle
\!\!\Big\langle [x,z]=[y,z]=0,\  [y,x] = \tau  z^2 \Big\rangle \!\!\Big\rangle
\right.
$$
The algebra $A^\tau$ has a natural grading defined by $\deg x = \deg y
= \deg z = 1$ and $\deg  \gamma = 0,$ for any $\gamma \in \Gamma$. Thus
$A^\tau  =  \bigoplus_{n \geq  0}  A^\tau_n$  is  a positively  graded
algebra, such that $A^\tau_0 = \C\Gamma$.
We will see in Appendix B that the algebra $A^\tau$
is strongly regular in the sense of Definition \ref{cond}.

\begin{example}\vi
If $\Gamma$ is trivial then $\tau$ reduces to a complex number. The
corresponding algebra, is a non-commutative deformation of the
polynomial algebra $\C[x,y,z]$.  

\vii If $\Gamma$ is
arbitrary and $\tau = 0$ then the algebra $A^\tau$ is the smash
product of the polynomial algebra $\C[x,y,z]$ with the group
$\Gamma$. 
\end{example}

We set $\PP^2_{\!_\G}=\proj(A^\tau),$ and write
$\coh(\PP^2_{\!_\G})=\qgr(A^\tau),$ for the corresponding category
of `coherent sheaves' on the non-commutative scheme $\PP^2_{\!_\G}$.
Further, let $\,\PP^1_{_\G}=\proj(\C[x,y]\#\Gamma)\,$
 be a noncommutative variety corresponding to
the graded algebra $\,\C[x,y]\#\Gamma$. The projection:
 $A^\tau\onto A^\tau/A^\tau z\cong\C[x,y]\#\Gamma$
 can be considered as a closed
embedding $i:\PP^1_{_\G}\into\PP^2_{\!_\G}$. Let
$i^*:\coh(\PP^2_{\!_\G})\to\coh(\PP^1_{_\G})$ and
$i_*:\coh(\PP^1_{_\G})\to\coh(\PP^2_{\!_\G})$ be the corresponding
pull-back and push-forward functors given in
terms of modules by the formulas
$$
i^*(\pi(M)) = \pi(M/Mz),\qquad
    i_*(\pi(N)) = \pi(N)\;,\;\text{ with $z$ acting by zero.}
$$

\begin{definition}
Let $\bcm_{_\G}^\tau(V,W)$ be the set of isomorphism classes of coherent
torsion free sheaves $E$ on $\PP^2_{\!_\G}$ with an isomorphism
$i^*E\cong W\gtimes\CO$ {\rm(}{\em framing of type $W$}{\rm)} and such
that $H^1(\PP^2_{\!_\G},E(-1))\cong V$,
as $\Gamma$-modules.
\end{definition}

In \S4 we will give a description of torsion free sheaves
on  $\PP^2_{\!_\G}$ in terms of monads, similar to the
well-known description of vector bundles on the commutative
$\PP^2$, cf. [OSS]. We will see that the linear algebra data
given by a monad  associated with
a sheaf on $\PP^2_{\!_\G}$ is nothing but a Nakajima quiver.
This way one obtains a bijection between the moduli space
$\bcm_{_\G}^\tau(V,W)$ and the  Nakajima quiver variety,
see Theorem \ref{quiv-sheaf} below.

\subsection{Quiver Varieties}
We keep the above setup, in particular,
$\Gamma$ is a finite subgroup of $SL_2(\C)$ and  $L$ is the  
tautological 2-dimensional representation of $\Gamma$. 
 Given a pair $(V, W)$ of finite dimensional
$\Gamma$-modules, and an element $\tau\in \ZZ(\CG)$,  
consider the locally closed subvariety of {\em quiver data}
\begin{equation}\label{bbm}
\bbm^\tau_{_\G}(V,W) \subset \Hom_{_\G}(V,V\otimes_\C L)
\;\bplus\;\Hom_{_\G}(W,V)\;\bplus\;\Hom_{_\G}(V,W)
\end{equation}
formed by all triples $(B,I,J)$ satisfying the following two
conditions: 
\begin{description}
\item[Moment Map Equation:] $[B,B] + IJ = \tau|_{_V}$\;\;\;
(here $[B,B]\in \End_{_\G}\!\!V\,\,\btimes\bigwedge^2\!L\simeq
\End_{_\G}V$)
\item[Stability Condition:] if $V'\subset V$ is a $\Gamma$-submodule 
such that $B(V')\subset V'\otimes L$ and 
\newline $I(W)\subset V'$ then $V'=V$.
\end{description}

The group $G_{_\G}(V)$ of $\Gamma$-equivariant
automorphisms of $V$ acts on
$\bbm^\tau_{_\G}(V,W)$ by the formula:
$
g(B,I,J) = (gBg^{-1},gI,Jg^{-1}).
$
Note that the $G_{_\G}(V)$-action is free,
due to the stability condition.

\begin{definition}
The quiver variety  is defined as the Geometric Invariant
Theory
quotient:
$\FM_{_\G}^\tau(V,W)=\bbm^\tau_{_\G}(V,W)/G_{_\G}(V)$. 
\end{definition}
Thus, the quiver variety above depends only on $[V]$,
the isomorphism class of the $\G$-module $V$, but abusing the notation 
we will write $\FM_{_\G}^\tau(V,W)$ rather than
$\FM_{_\G}^\tau([V],W)$.

The relation of the above definition of quiver variety with the
original definition of Nakajima is provided by  McKay
correspondence.  Recall  that  McKay
correspondence, cf. \cite{Re}, establishes a bijection between affine
Dynkin graphs of $ADE$-type,  and finite subgroups  $\Gamma \subset
SL_2(\C)$
  (up to
conjugacy). Given  such a  subgroup $\Gamma$, the  corresponding affine
Dynkin graph can be recovered as follows.

Let  $R_i$,  $i  = 0,  \ldots,  n$  be  the  complete set  of  the
isomorphism classes of complex
irreducible representations of $\Gamma$.
 For  any $i,  j \in  \{0, \ldots, n\},$  let $a_{ij}$  be the
multiplicity  of $R_j$  in the  $\Gamma$-module $R_i  \otimes_{\C} L$,
where $L$ is the tautological 2-dimensional representation.
Since $L$ is  self-dual, we have $a_{ij} =  a_{ji}$.  Following McKay,
we attach to  $\Gamma$ a graph $Q=Q(\Gamma)$  by taking
$\{0,\ldots, n\}$ as  the set of vertices of $Q$, and  connecting $i$ and $j$
by $a_{ij}$ edges.    Then  $Q(\Gamma)$   is  an affine   Dynkin  graph
of ADE-type,
with Cartan matrix $a(Q)=\|a_{ij}\|$.
The Nakajima's quiver variety corresponding to $Q(\Gamma)$
 coincides with the variety $\FM_{_\G}^\tau(V,W)$ 
defined above. 

\begin{remark} The case: $\Gamma = \{1\}$ is somewhat degenerate.
In this case $Q(\Gamma)$ should stand for the quiver
with one vertex, and {\it one} (rather than two) loop.
$\quad\lozenge$
\end{remark} 

Our first important result relates torsion free sheaves on 
$\PP^2_{\!_\G}$ with the Nakajima quiver variety:

\begin{theorem}\label{quiv-sheaf}
There exists a natural bijection:
$\FM^\tau_{_\G}(V,W)\iso \bcm_{_\G}^\tau(V,W)$.
\end{theorem}

\begin{remark}
It is natural to expect that the bijection of the
theorem establishes an isomorphism of algebraic varieties.
We are unaware, however, of any construction of a scheme (or stack)
structure
on the set $\bcm_{_\G}^\tau(V,W)$,
since the
formalism of moduli spaces of sheaves on noncommutative varieties is
not yet developed.$\quad\lozenge$
\end{remark}

\nc{\BQ}{{\bar{Q}}}
\nc{\CBQ}{{\C\BQ}}
\medskip

 The  above theorem  can be viewed  in the following  context. 
 The Kleinian singularity, 
  $\C^2/\Gamma,$ has a semiuniversal deformation, a family 
$\{\X_\tau\}$ of surfaces parametrised by the
points $\tau$ in the 
vector space ${\mathfrak{h}} =
\ZZ^0(\CG) \subset
\ZZ(\CG)$, the codimension one
 hyperplane in $\ZZ(\CG)$
formed by all central elements which have trace zero in
the regular  representation $\CG$. Thus, $\X_0=\C^2/\Gamma$.
Now, let $W=\triv$ be the
trivial $\Gamma$-module, and  $V=\CG$, the regular
representation. For each $\tau
\in {\mathfrak{h}}$ put $\widetilde{\X}_\tau=\FM_{_\G}^\tau(\CG,\triv)$.
According to Kronheimer [Kr],
the family $\,\{\widetilde{\X}_\tau\}_{\tau
\in \mathfrak{h}}\,$ gives a simultanious
resolution: $\widetilde{\X}_\tau \onto \X_\tau$
of the semiuniversal deformation, $\{\X_\tau\}_{\tau
\in \mathfrak{h}}$. In particular,
$\widetilde{\X}_0\onto\C^2/\Gamma$ is the minimal resolution.

Kronheimer and Nakajima  [KN] showed further that,
for arbitrary $V$ and $W$, and $\tau\in \h=\ZZ^0(\CG)$,
the  quiver variety
$\FM_{_\G}^\tau(V,W)$ 
 is isomorphic to the moduli space of
anti-self-dual connections on 
$\widetilde{\X}_{-\tau}$ with a suitable decay
condition at infinity. If $\tau=0$, then
there is an alternative, purely algebraic, interpretation
of $\FM_{_\G}^\tau(V,W)$ as the moduli space of
framed torsion-free sheaves (satisfying an appropriate
stability condition)  on 
$\widehat{\X}_0$, a natural projective completion of
$\widetilde{\X}_0$.
Such an  interpretation of $\FM_{_\G}^\tau(V,W)$
in terms of  moduli spaces of
framed torsion-free sheaves on $\widehat{\X}_\tau$
doesn't seem to be possible, however, for other $\tau\in \h
\smallsetminus \{0\}$.

A way to circumvent this difficulty has been suggested by
 Kapranov and  Vasserot, who
  proved  in \cite{KV}  that considering framed 
 torsion-free sheaves  on
$\widehat{\X}_0$ is equivalent
  to considering $\Gamma$-equivariant framed sheaves on $\PP^2$.
  Furthermore, Varagnolo and   Vasserot proved in [VV]
that $\Gamma$-equivariant framed torsion-free
sheaves on  $\PP^2$ (equivalently, framed torsion-free
sheaves on  the non-commutative
  $\proj$-scheme
 corresponding to the graded  algebra $A^0=\C[x, y, z] \# \Gamma$)
are parametrised by the points of $\FM_{_\G}^0(V,W)$.
Thus, our theorem above is a very natural extension of \cite{VV}.

Notice further, that while the spaces $\widetilde{\X}_\tau$ are only
defined for $\tau\in \h=\ZZ^0(\CG)$, 
the  quiver  varieties $\FM_{_\G}^\tau(V,W)$
 exist \textit{for  all} $\tau
  \in \ZZ(\CG)$.
This discrepancy has  been resolved  by Crawley-Boevey  and  Holland,  who
 constructed in [CBH] a family of associative algebras  $\A_\tau$
parametrised
by all $\tau
 \in  \ZZ(\CG)$, and such that the algebra $\A_\tau$
is isomorphic to the coordinate ring of
the variety $\mathcal{\X}_\tau$, if   $\tau  \in  \mathfrak{h}$,  
and is \textit{non-commutative} otherwise. 
From this perspective, 
our theorem can also be viewed as describing framed torsion-free
sheaves on a projective completion of $\Spec\A_\tau$,
the non-commutative affine scheme corresponding to
the algebra  $\A_\tau$.
\medskip

Next we introduce the algebra
$$
\B^\tau = A^\tau/(z-1)\cd A^\tau\;\simeq\;
\bigl(\C\langle x,y\rangle \# \Gamma\bigr)
\big/\langle\!\langle [y,x] = \tau \rangle\!\rangle,
$$
which can be thought of as the `coordinate ring' of 
$\PP^2_{\!_\G}\smallsetminus \PP^1_{_\G}$, 
a non-commutative affine plane.

The algebra $\B^\tau$ comes equipped with
a natural increasing filtration, $F_\bullet\B^\tau$,
such that $F_0\B^\tau=\CG$. It follows from [Q] that
the assignment: $R\mapsto R\otimes_{_\CG}\B^\tau$
gives a natural isomorphism:
$K(\G)\iso K(\B^\tau),$ between
the Grothendieck group of the category of 
finite dimensional $\Gamma$-modules and that of
finitely generated projective
$\B^\tau$-modules. Further,
the map assigning the integer
$\dim R$ to a finite dimensional $\Gamma$-module $R$
extends to a canonical group homomorphism $\dim: K(\G)\to \Z$.
Composing this with the isomorphism $K(\B^\tau)\simeq K(\G)$,
one gets a function ${\widetilde{\dim}}: K(\B^\tau)\to \Z$.

The second major result of this paper is a classification
of isomorphism classes of finitely generated
projective $\B^\tau$-modules $N$
such that ${\widetilde{\dim}} N=1$. In order to formulate it
we need the following technical result whose proof is given in \S8.

\begin{proposition}\label{easy} Let $R\in K(\G)$
be such that $\dim R=1$. Then there exist
uniquely determined (up to isomorphism)
$\G$-modules $W$ and $V$, such that in $K(\G)$ we have:
$R= [W]+ [V]\cdot([L]-2[\triv])$ and, moreover,
 $\dim W=1$ and  $V$ does not
contain the regular representation as a submodule.
\end{proposition}

The
Theorem below
is a generalization of  a conjecture of
Crawley-Boevey and Holland, see~\cite{BLB}, Example 5.7.
The conjecture corresponds, when reformulated in terms
of deformed preprojective algebras, see \S6,
to the special case of the Theorem: $R=[W]\,,\,V=~0.$

\begin{theorem}\label{main_intro}
Assume that $\tau\in \ZZ(\CG)$  is generic in the sense
of [CBH], i.e., it does not belong to root hyperplanes.
Let $R\in K(\G)$ be such that $\dim R=1$, and
$V\,,\,W$ the $\G$-modules attached to this class by Proposition
\ref{easy}. Then there exists
a natural bijection:
$$
\left\lbrace
\begin{array}{l}\mbox{Isomorphism classes of finitely generated 
projective}\\
\mbox{$\B^\tau$-modules $N$ such that 
$[N]=R$ in $K(\B^\tau)=K(\Gamma)$}
\end{array}
\right\rbrace
\;\;\simeq\;\;
\bigsqcup_{k=0}^\infty\; \FM_{_\G}^\tau(V\oplus\CG^{\oplus k}\,,\,W)
\,.
$$
\end{theorem}

Recall that to define "root hyperplanes" in $\ZZ(\CG)$,
one has to identify this space with the underlying vector space
of the affine root system, $\Delta_{aff}\subset\ZZ(\CG)^*$,
associated to the Dynkin graph $Q(\G),$
see also Sections 5 and 6 of \cite{CBH}.
The vertices of $Q$ may be identified with the base of the space $\ZZ(\CG)^*$
formed by the irreducible characters of $\G$,
 and the root system $\Delta_{aff}$
is given in this basis by the Cartan
 matrix $a(Q)=\|a_{ij}\|$ of the graph $Q$.\medskip

To prove Theorem \ref{main_intro} we observe first that,
 for generic $\tau\in \ZZ(\CG)$,
any torsion-free sheaf on $\PP^2_{\!_\G}$ is in effect 
locally free. We show further that `restriction' from
$\PP^2_{\!_\G}$ to $\PP^2_{\!_\G}\smallsetminus \PP^1_{_\G}$, 
gives, for generic $\tau\in \ZZ(\CG)$, a bijection
between locally free sheaves on $\PP^2_{\!_\G}$ and
$\PP^2_{\!_\G}\smallsetminus \PP^1_{_\G}$, respectively.
The latter objects being nothing but finitely generated projective
$\B^\tau$-modules, the result follows from Theorem \ref{quiv-sheaf}.
\medskip

Let $\Gamma=\{1\}$ so that: $R=\triv$, and  $\ZZ(\CG)=\C$.
The parameter $\tau\in \C$ 
 is generic in the above sense if and only if
$\tau\neq0$. The algebra $\B^\tau$ becomes,
for $\tau\neq0$, the first Weyl algebra with two generators,
$x,y$, subject to the relation: $[y,x]=\tau\cdot 1$.
Further, the variety $\FM_{_\G}^\tau(\CG^{\oplus k}, \triv)$
becomes in this case the Calogero-Moser variety,
$\FM^\tau(\C^k)=\FM_{_{\G=\{1\}}}^\tau(\C^{\oplus k}, \C)$,
introduced in [KKS] and studied further in [W]. Explicitly,
we have:
$$\FM^\tau(\C^k)\;=\;\big\lbrace{(B_1,B_2)\in \gln\times\gln
\;\;\Big|\;\;
[B_1,B_2]-\tau
\cdot\id_{_{\C^n}}
=\;\text{rank } 1 \text{ matrix}\big\rbrace}/\text{Ad}\,GL_n
\,.
$$

Since any rank 1 projective module over the Weyl algebra  $\B^\tau$
is isomorphic
to a right ideal in  $\B^\tau$,
Theorem \ref{main_intro} reduces in the
special case of the trivial group $\G$ to the following
result proved recently by Berest-Wilson [BW].

\begin{corollary}\label{calogero}
The set of isomorphism classes of right ideals
(viewed as $\B^\tau$-modules) in 
the Weyl algebra
$\B^\tau$ is in a natural bijection 
with the set: 
$\;
\bigsqcup_{k=0}^\infty\; \FM^\tau(\C^k).
$
\end{corollary}

The proof of this result given in [BW]
is totally different from ours, and relies heavily 
on the earlier results of [W] and [CaH].

\section{Sheaves on a noncommutative {\sf{Proj}}-scheme}
\setcounter{equation}{0}

In this section we fix a strongly regular graded
$\C$-algebra 
 $A=\bigoplus_{k\geq 0}\,A_k$, set
$X=\proj A$, and write $\coh(X)=\qgr(A)$.
Recall that the quotient category $\qgr(A)$ is defined as follows.
The objects of $\qgr(A)$ are the same as those of $\gr(A)$, while
$$
\Hom_{\qgr(A)}(M,N) = \lim\limits_{\longrightarrow}\Hom_{\gr(A)}(M',N),
$$
where  the limit is taken  over all submodules  $M'\subset M$, with
$M/M'\in\tor(A)$. 

Recall the integer $d$ entering the definition \ref{reg}
of a strongly regular
algebra and the Serre duality property (\ref{serre})

\begin{proposition}
For any coherent sheaves $E$ and $F$ on $X$ we have
$$
\Ext^{>d}(E,F) = 0,\qquad \CExt^{>d}(E,\CO)=0.
$$
\end{proposition}
{\sl Proof:} 
The first part follows from Serre duality and the second follows
from the first and from (\ref{cextp}).
\qed\medskip

\begin{proposition}\label{resolution}
Any coherent sheaf $E$ on $X$  admits a resolution of the form
\begin{equation}\label{oresol}
\dots \to \CO(-n_k)^{\oplus m_k} \stackrel{\phi_k}\lra \dots 
\stackrel{\phi_1}\lra \CO(-n_0)^{\oplus m_0} \to  E \to 0.
\end{equation}
\end{proposition}
{\sl Proof:} 
Follows immidiatelty from the ampleness property.
\qed\medskip

Recall 
  the notions of 
locally free and torsion free sheaves on $X$, see
Definition \ref{free}. We write $E^*$ instead of
$\CHom(E,\CO)$, for brevity.

\begin{proposition}\label{lfprop}
Let $\CE$, $\CE_1,\dots,\CE_d$ be locally free sheaves on $X$.

\newi1 For any complex:
$\CE_1\stackrel{\phi_1}\too\dots\stackrel{\phi_{d-1}}\too\CE_d,$
which is exact at the terms: $\CE_2$, \dots, $\CE_{d-1},$ 
the sheaf $\Ker\phi_1$ is locally free.

\newi2 The sheaf $\CE^*$ is locally free.

\newi3 The canonical morphism: $\CE\to\CE^{**}$ is an isomorphism.

\newi4 For any $k$, the sheaf $\CE$ has a resolution of the form 
\begin{equation}\label{rresol}
0 \to \CE \to \CO(n_0)^{\oplus m_0} \to \dots 
\to \CO(n_k)^{\oplus m_k} \to  \dots
\end{equation}

\newi5 If $d\le 2$, then for any sheaf $E$,
 the sheaf $E^*$ is locally free.\footnote{we will usually denote
arbitrary sheaves by roman letters, and locally free sheaves by
script letters.} 
\end{proposition}
{\sl Proof:}
\newi1 Let $K=\Ker\phi_1$ and $C=\Coker\phi_{d-1}$. Applying
the functor $\CHom(-,\CO)$ to the sequence
$0\to K\to \CE_1\to\dots\to\CE_d \to C\to 0$ we get a spectral sequence
converging to zero. The first term of the sequence looks as
$$
\xymatrix@R-15pt@!C{
{\CExt}^d(K,\CO) & 0             & {\dots}      &     0          & 
{\CExt}^d(C,\CO) \ar@{-->}[ddllll]\\
{\vdots}       & {\vdots}      & {\vdots}     & {\vdots}       & {\vdots} \\
K^*            & \CE_1^* \ar[l]& {\dots}\ar[l]& \CE_d^* \ar[l] & C^* \ar[l]
}
$$
Hence $\CExt^{>0}(K,\CO)=0$.

\newi2 Choose a resolution (\ref{oresol}) for $\CE$ with $k=d$ and let
$\CF=\Ker\phi_d$. Then by (1) the sheaf $\CF$ is locally
free. Applying the functor $\CHom(-,\CO)$ to this resolution we get a
resolution 
$$
0 \to \CE^* \to \CO(n_0)^{\oplus m_0} \to \dots 
\to \CO(n_d)^{\oplus m_d} \to  \CF^* \to 0
$$
Hence $\CE^*$ is locally free by (1).

\newi3 For  $\CE$, choose a resolution of type
(\ref{oresol})  with $k=d$, and
apply  to it the functor $\CHom(\CHom(-,\CO),\CO)$. Since the 
canonical morphism:
$\CO(n)\to\CO(n)^{**}$ is an isomorphism it follows that:
$\CE\iso\CE^{**}$ is also an isomorphism. 

\newi4 For  $\CE^*$,
choose a resolution  of type (\ref{oresol}), and apply
$\CHom(-,\CO)$.

\newi5 For  $\CE$, choose a resolution  of type 
(\ref{oresol})  with $k=1$. Apply to it
the functor $\Hom(-,\CO)$, and use (1).~\qed\medskip

\begin{proposition}\label{tfprop}
If $E$ is a torsion free sheaf then $\CExt^d(E,\CO) = 0$.
\end{proposition}
{\sl Proof:}
Given $E$, choose an embedding: $E\into\CE$  into a locally free sheaf.
Applying the functor $\CHom(-,\CO)$ to the exact sequence
$0\to E\to\CE\to\CE/E\to0$ we get an epimorphism:
$\CExt^d(\CE,\CO)\onto\CExt^d(E,\CO)\to0$.
\qed\medskip


\begin{lemma}\label{exti}
For $n\gg0$, there is a canonical isomorphism:
$$
H^0\Bigl(X,\bigl(\CExt^k(E,\CO)\bigr)(n-l)\Bigr)\cong
H^{d-k}(X,E(-n))^\vee
\quad,\quad\forall k\geq 0.
$$
\end{lemma}
{\sl Proof:}
It is clear that
$H^p\Bigl(X,\bigl(\CExt^k(E,\CO)\bigr)(n-l)\Bigr)
=H^p\Bigl(X\,,\,\CExt^q\bigl(E(-n)\,,\,\CO(-l)\bigr)\Bigr)$.
Further, there is a standard spectral sequence for the
derived functor of composition:
$$
E_2^{p,q}=H^p\Bigl(X,\CExt^q\bigl(E(-n),\CO(-l)\bigr)\Bigr) \Longrightarrow
\Ext^{p+q}\bigl(E(-n),\CO(-l)\bigr).
$$
On the other hand, for $n\gg0$ we have:
$H^{>0}\Bigl(X,\bigl(\CExt^k(E,\CO)\bigr)(n-l)\Bigr)=0$,
hence 
the spectral sequence degenerates and
$H^0\Bigl(X,\bigl(\CExt^k(E,\CO)\bigr)(n-l)\Bigr)=\Ext^k\bigl(E(-n),
\CO(-l)\bigr)$.
Finally, Serre duality gives
$$
\Ext^k\bigl(E(-n),\CO(-l)\bigr)\;\cong\; H^{d-k}(X,E(-n))^\vee.\qquad\square
$$

\begin{definition}
A sheaf $F$ on $X$ is called an {\em Artin sheaf} if
$H^{>0}(X,F(k))=0\,,\,\forall k\in\Z$.
\end{definition}
\begin{remark} 
In the situations studied in this
paper,
any Artin sheaf has finite length. We don't know if this is true in
general, as well as whether  a subsheaf of an Artin sheaf is Artin.
$\quad\lozenge$
\end{remark}

\begin{proposition}\label{artin} Assume that $d = \dim X\le 2$. 

\newi1 If $F$ is an Artin sheaf and $\CE$ is locally free
then $\Ext^{>0}(\CE,F)=0$.

\newi2 For any sheaf $E$ the sheaf $\CExt^d(E,\CO)$ is Artin.

\newi3 $F$ is an Artin sheaf if and only if $\CExt^{<d}(F,\CO)=0$.

\newi4 If $F$ is an Artin sheaf then $\CExt^d(F,\CO)=0$ implies $F=0$.

%
\newi5 Any extension of Artin sheaves is Artin.
\end{proposition}

\noindent
{\sl Proof:}
\newi1 First note that $\Ext^p(\CO(k),F)=H^p(X,F(-k))=0$, for $p>0$.
Now, assume $\CE$ is locally free and take a resolution
of $\CE$ of type (\ref{rresol}) with $k=d$.
Applying the functor $\Hom(-,F)$ to this resolution we 
deduce that $\Ext^{>0}(\CE,F)=0$.

\newi2 
We have a spectral sequence
$$
E_2^{p,q} = H^p\bigl(X\,,\,
\CExt^q(E,\CO(k))\bigr) \Longrightarrow \Ext^{p+q}(E,\CO(k)).
$$
It is easy to see that any nonzero element in
$H^p\bigl((X\,,\,\CExt^d(E,\CO(k))\bigr)$
 with $p>0$ gives a nonzero contribution
in $\Ext^{d+p}(E,\CO(k))$. But the latter space is zero for $p>0$.
Hence $H^p\bigl((X\,,\,\CExt^d(E,\CO(k))\bigr)=H^p\bigl((X\,,\,
\CExt^d(E,\CO)(k)\bigr)=0$,
that is $\CExt^d(E,\CO)$ is an Artin sheaf.

\newi3 If $F$ is an Artin sheaf then it follows from
Lemma \ref{exti} and
ampleness that $\CExt^{<d}(F,\CO)=0$. 
Conversely, assume that $\CExt^{<d}(F,\CO)=0$.
We have a spectral sequence
$$
E_2^{p,q} = \CExt^p(\CExt^{-q}(F,\CO),\CO) \Longrightarrow E_\infty^i =
\begin{case} F, &\text{if $i=0$}\\ 0, &\text{otherwise}\end{case}
$$
Since $\CExt^{<d}(F,\CO)=0,$ the spectral sequence degenerates at the
second term and yields an isomorphism:
 $F\cong\CExt^d(\CExt^d(F,\CO),\CO)$.
Hence $F$ is an Artin sheaf by (2).

\newi4 It follows from the proof of (3) that for any Artin sheaf $F$
we have a canonical isomorphism:
$$
F\cong\CExt^d(\CExt^d(F,\CO),\CO).
$$
Hence, $\CExt^d(F,\CO)=0$ implies $F=0$.\quad (5) is clear.
\quad$\square$


\section{Torsion-free sheaves on $\PP^2_{\!_\G}$}
\setcounter{equation}{0}

In this section we set $A^\tau
=  \bigl(\C  \langle  x,y,z  \rangle  \#  \Gamma\bigr)  \big/  \big\langle
\!\big\langle [x,z]=[y,z]=0,\  [y,x] = \tau  z^2 \big\rangle \!\big\rangle
,$ see Definition \ref{Atau}. We write:
$A^\tau
= \bplus_{\! i\geq 0}\;A_i,\,$ and put $\,\PP^2_{\!_\G}=\proj(A^\tau)$.
\subsection{Beilinson Spectral Sequence}
The noncommutative projective plane $\PP^2_{\!_\G}$ shares a lot of
properties with the commutative projective plane. In particular, we
have,
see Appendix B for the definition of ${}^*  \!  A_j$:
\begin{equation}\label{cohoi}
H^p(\PP^2_{\!_\G},\CO(i))=\begin{case} 
A_i, & \text{if $p=0$ and $i\ge0$}\\
{}^*  \!  A_{-i-3},  &   \text{if  $p=2$  and  $i\le-3$}\\  0\quad,  &
\text{otherwise}
\end{case}
\end{equation}

Our approach to the classification of torsion free sheaves on
$\PP^2_{\!_\G}$ mimics the standard approach,
see [OSS],[KKO], to the study of coherent sheaves on
$\PP^2$ by means of Beilinson's spectral sequence,
see Appendix B. The latter allows one to describe coherent
sheaves in terms of linear algebra data, sometimes called
the `ADHM-equations'.  

Beilinson's spectral sequence for a general Koszul algebra, given in
Appendix B, simplifies considerably in the `2-dimensional' case of
$\PP^2_{\!_\G}$. Specifically,
let $T$ denote the coherent sheaf on $\PP^2_{\!_\G}$  defined
by either of the following two exact sequences:
\begin{equation}\label{deft} 
{}^!\!A_0\gtimes\CO \into {}^!\!A_1\gtimes\CO(1) \onto T \qquad,\qquad
T \into {}^!\!A_2\gtimes\CO(2) \onto {}^!\!A_3\gtimes\CO(3), 
\end{equation}
where ${}^!\!A_k$ stand for the graded 
components of the Koszul dual algebra,
see Appendix B.
Then, for any sheaf $E$, there is a spectral sequence with the $
E_1^{p,q}$-term looking like a 3-term complex:
$$
E_1^{p,q} = \Big\{
\Ext^q(\CO(1),E)\gtimes\CO(-2) \to
\Ext^q(T(-1),E)\gtimes\CO(-1) \to
\Ext^q(\CO,E)\gtimes\CO
\Big\}
$$
The above complex corresponds to the groups
$
E_1^{p,q}$ with $p=-2,-1,0$, all other groups being zero. 
This spectral sequence converges to
(see e.g., \cite{OSS} for details in the commutative case):
$$
E_1^{p,q} \Longrightarrow E_\infty^{p+q} =
\begin{case} E, & \text{for $p+q=0$}\\ 0, & \text{otherwise}\end{case}\,.
$$

Recall that $K(-)$ stands for the Grothendieck
group of an abelian category. 

\begin{corollary}\label{k0}
We have $K(\coh(\PP^2_{\!_\G})) \cong K(\Gamma)^{\oplus 3}$. In
particular, $K(\coh(\PP^2_{\!_\G}))$ is a free $\Z$-module with a
basis given by classes of sheaves $R\gtimes\CO(k)$, where $R$
runs through the set of isomorphism classes of irreducible
$\Gamma$-modules and $k\in\{-2,-1,0\}$. $\quad\square$
\end{corollary}

For any $E\in\coh(\PP^2_{\!_\G})$ we define the {\em Hilbert function}
of $E$ as
$$
h_E(t) = \sum\nolimits_{p=0}^2 \;\,(-1)^p\cdot
\dim_\C H^p(\PP^2_{\!_\G},E(t)).
$$

\begin{lemma}\label{hilbert}
For any $E\in\coh(\PP^2_{\!_\G})$, the function $h_E(t)$ is a polynomial
 in $t$ (of degree $\le 2$)  of the form
$$
h_E(t) = r\mbox{$\frac{t^2}2$} + \dots\quad,\quad r=1,2,3,\ldots\;\,.
$$
\end{lemma}
{\sl Proof:}
It is clear that $h_E(t)$ depends only on the class of $E$ in
$K(\coh(\PP^2_{\!_\G}))$. Thus by Corollary \ref{k0} it suffices to compute the
Hilbert function only for the sheaves of the form
$R\gtimes\CO(p)\,,\,p\in\Z$, where $R$ is a $\Gamma$-module.
Further,  from~(\ref{cohoi}) one deduces  that for $t+k\ge 0$ we have:
$H^{>0}(\PP^2_{\!_\G},R\gtimes\CO(t+k)) = 0$. Therefore, we find:
$$
\begin{array}{l}
h_{R\gtimes\CO(k)}(t) =  
\dim_{_\C} H^0(\PP^2_{\!_\G}\,,\,R\gtimes\CO(t+k)) =
\medskip\\\qquad\qquad =
\dim_{_\C} \bigl(R\gtimes A^\tau_{t+k}\bigr) =
\dim_{_\C} \Bigl(R\gtimes\bigl(\CG\otimes_{_\C} Sym^{t+k}\langle x,y,z
\rangle\bigr)\Bigr)= 
\medskip\\\qquad\qquad\displaystyle =
\dim_{_\C}\bigl( R\otimes_{_\C} Sym^{t+k}\langle x,y,z\rangle\bigr) =
\frac{(t+k+1)(t+k+2)}2\cdot \dim_{_\C} R
\end{array}
$$
A similar computation shows that this formula  also holds  for $t+k<0$.
\qed\medskip

\begin{definition}
Define the rank $r(E)$ of  $E\in\coh(\PP^2_{\!_\G})$ to be the
leading coefficient $r$ of the Hilbert polinomial 
$h_E(t)$, see Lemma {\rm\ref{hilbert}}.
\end{definition}

\begin{corollary}\label{rri}
The rank, $r(-)$, is a well-defined linear function on
$K(\coh(\PP^2_{\!_\G}))$. Moreover, 
$r(R\gtimes\CO(i)) = \dim_{_\C} R$, for any $\Gamma$-module $R$ and
any $i\in\Z$.$\quad\square$
\end{corollary}

\subsection{Sheaves on $\PP^1_{_\G}$}

Recall that in \S1.2 we have defined the noncommutative projective 
line as the  $\proj$-scheme corresponding to the algebra
$\C[x,y]\#\Gamma$. The following result is clear

\begin{proposition}\label{v1}
The embedding of graded algebras $\C[x,y]\subset\C[x,y]\#\Gamma$
induces an equivalence of categories
$$
\coh(\PP^1_{_\G}) \to \coh\nolimits_{_\G}(\PP^1),
$$
where $\coh_{_\G}(\PP^1)$ is the category of $\Gamma$-equivariant
coherent sheaves on the commutative $\PP^1$.\qed
\end{proposition}

Note that the  equivalence of Proposition \ref{v1}
 commutes with the functors
$\Ext^p$, $\CExt^p,$ and $H^p$. It follows that a sheaf $E$
on $\PP^1_{_\G}$ is locally free (resp. Artin) if and only if  it is locally free
(resp. Artin) as a $\Gamma$-equivariant sheaf on $\PP^1$. For any
sheaf $E\in\coh(\PP^1_{_\G}),$ we denote by $h_E(t)$, resp. $r(E)$, the
Hilbert polynomial, resp.   the rank, of $E$, considered as a
$\Gamma$-equivariant sheaf on $\PP^1$. 

\begin{corollary}\label{plg}
\newi1 For any coherent sheaf $E$ on $\PP^1_{_\G}$ we have
$E\cong F\oplus\CE$, where $F$ is an Artin sheaf
and $\CE$ is a locally free sheaf.

\newi2  Any locally free sheaf on $\PP^1_{_\G}$ has the form:
$\CE=\bplus_{k}\, \Bigl(R_k\gtimes\CO(k)\Bigr)$,
 for certain $\Gamma$-modules~$R_k$.

\newi3 If $\CE$ is a locally free sheaf on $\PP^1_{_\G}$ and
$r(\CE)=1$, then $\CE\cong R\gtimes\CO(k),$ 
where $R$ is a 1-dimensional $\G$-module.$\quad\square$
\end{corollary}

\subsection{The functors $i^*$ and $i_*$}

Recall that in  section 1.2 we have defined the functors
$$
i_*:\coh(\PP^1_{_\G})\to\coh(\PP^2_{\!_\G})\quad\text{and}\quad
i^*:\coh(\PP^2_{\!_\G})\to\coh(\PP^1_{_\G}).
$$ 
It is clear that $i_*$ is the
right adjoint of the
functor $i^*$. 
It follows from the definition that 
$i^*$ is right exact. We denote by $L^pi^*$ the left derived
functor.


\begin{proposition}\label{iprop}

\newi1 The functor $i_*$ is exact and faithful.

\newi2 For any coherent sheaf $E$ on $\ppgt\tau$ 
there is a canonical  exact sequence:
$$
0\too 
i_*L^1i^*E\too  E(-1)\stackrel{z\cdot}\too E \too i_*L^0i^*E\too 0
$$

\newi3 For any $E$ we have $L^{>1}i^*E=0$ and the functor $L^1i^*$ is
left exact.

\newi4 If $\CE$ is locally free then $L^{>0}i^*\CE=0$.

\newi5 If $E$ is torsion free then $L^{>0}i^*E=0$.

\newi6 If $\CE$ is a locally free sheaf on $\ppgt\tau$ then $i^*\CE$
is locally free.

\newi7 If $E$ is torsion free then $r(i^*E) = r(E)$.

\newi8 For any sheaf $E$ on $\ppg$,
the  adjunction morphism: $E\to i_*i^*E$ is an epimorphism.

\newi9 We have: $i^*i_*E=E$, $L^1i^*i_*E=E(-1),$ for any sheaf
$E$ on $\plg$.

\newi{10} If $i^*E=0$ then $E\cong E(1)$, and $E$ is an Artin sheaf.
\end{proposition}
{\sl Proof:}
\newi1 This claim becomes clear when translated into the module language.

\newi2 Since $i_*$ is exact it suffices to check that:
$i_*i^*E\cong E\otimes i_*\CO_{\plg}$ (which is clear
from the point of view of modules), and to apply 
the resolution
$$
0 \to \CO(-1) \stackrel{z}\lra \CO \to i_*\CO_{\PP^1_{_\G}} \to 0. 
$$

\newi3 Follows from (2) and (1).

\newi4 It is clear that $L^{>0}i^*\CO(n)=0$.
Choosing an embedding $\CE\to\CO(n)^{\oplus m}$ like in (\ref{rresol})
and applying (3) we obtain the claim.

\newi5 Given a torsion free sheaf $E$,
choose an embedding: $E\into\CE,$
 into a locally free sheaf.
Now, apply (3) and (4).

\newi6 Note that derived functor $L^\cdot  i^*$
commutes with the derived functor
$\CExt^\cdot(-,\CO)$. Applying (4) to the  locally free sheaf 
$\CE$ we
conclude that $\CExt^{>0}(i^*\CE,\CO_{\plg})=0$.
Hence $i^*\CE$ is
locally free.

\newi7 Since $L^{>0}i^*E=0$ by (5) it follows that:
$\;
h_{i^*E}(t) = h_E(t) - h_E(t-1) = r(E)t + \dots
\,.$ 
The claim follows.

\newi8 Follows from (2).\quad (9)$\;$ Note that by (2) we have
$$
i_*i^*i_*E    =  \Coker (i_*E(-1) \stackrel z\to i_*E)
\qquad\mbox{and}\qquad
i_*L^1i^*i_*E =    \Ker (i_*E(-1) \stackrel z\to i_*E).
$$
On the other hand $z$ vanishes on $\plg$, hence
the morphism: $i_*E(-1) \stackrel z\to i_*E$ vanishes.
Thus $i_*i^*i_*E=i_*E$ and $i_*L^1i^*i_*E=i_*E(-1)$.
Hence by (1) we have: $i^*i_*E=E$, $L^1i^*i_*E=E(-1)$.

\newi{10} The equation: $i^*E=0$ implies by (2) that
 the multiplication by $z$ homomorphism: $E
\stackrel{z\cdot}{\too} E(1)$ is surjective. 
Let  $M$ be the  graded $A$-module corresponding to $E$.
 It
follows that,  for $k\gg0$, the
$z$-multiplication:
 $M_k\stackrel{z\cdot}{\too} M_{k+1}$ 
is a surjection, by the  ampleness. Since
$M$ is finitely generated it follows 
that $M_k\cong M_{k+1}$ for
$k\gg0$. Hence $E\cong E(1)$. But this implies that for any $k$ we
have $H^{>0}(\PP^2_{\!_\G},E(k))=0$. Thus, $E$ is an Artin sheaf.
\qed\medskip

\begin{definition}
A sheaf $E$ on $\ppgt\tau$ is said to be {\em supported on $\plg$}
if it admits an increasing finite
 filtration: $0 = E_0\subset E_1\subset\dots\subset E_n=E\,,\, E_i\in 
\coh(\ppgt\tau),$
such that $E_k/E_{k-1}=i_*F_k$, for some $F_k\in\coh(\plg)$.
\end{definition}

\begin{proposition}\label{supi}
\newi1 If $E$ is supported on $\plg$ then $r(E)=0$.

\newi2 Assume that $E$ is supported on $\plg$.
Let $E^0=E$ and, for $k>0$, 
 define $E^k$  inductively by
$\,E^{k+1} := \Ker(E^k\to i_*i^*E^k)$. Then $E^n=0$, for  $n\gg 0$.

\newi3 If $E$ is supported on $\plg$ and $i^*E$ is Artin
then $E$ is Artin.
\end{proposition}
{\sl Proof:}
\newi1 It is clear that $r(E)=\sum r(E_k/E_{k-1}) = \sum r(i_*F_k)$.
On the other hand, we have
$$
h_{i_*F}(t) = h_F(t) = r(F)\cdot t + \text{const},
$$
for any sheaf $F$ on $\plg$. Hence $r(i_*F)=0$.

\newi2 It suffices to prove the following claim:
{\it  Assume that
we have an exact sequence
$$
0 \to E' \to E \to i_*F \to 0.
$$
Then $\Ker(E\to i_*i^*E)$ is a subsheaf in $E'$.}
Indeed, if the claim is true then it is easy
to prove by induction that $E^k\subset E_{n-k}$,
where $E_\bullet$ is a filtration from the definition
of a sheaf supported on $\plg$. Hence $E_n\subset E_0 =0$,
and (2) follows.

We now  prove the claim. Since the canonical map:
$i_*F \to i_*i^*i_*F$ is an isomorphism by~\ref{iprop}~(9)
it follows that the projection: $E\onto i_*F$ factors through
the morphism: $E\to i_*i^*E$. Thus, the claim follows from the inclusion:
$$
\Ker(E\to i_*i^*E)\;\subset\;\Ker(E\to i_*F)\, =\, E'.
$$

\newi3 Assume that $E$ is supported on $\plg$ and that $i^*E$
is an Artin sheaf. We define inductively
 the sequence of sheaves $E^0=E\,,\,$
$E^{k+1}=\Ker(E^k\to i_*i^*E^k)\,,\,k=0,1,2,\ldots,$ as in (2). 
According to part (2), there exists $n>0$ such that
$E^n=0$. We prove by descending induction on $k$,
starting at $k=n$,
that $E^k$ is an Artin sheaf. Applying the functor $i^*$ to
the exact sequence
$$
0 \too E^{k+1} \too E^k \too i_*i^*E^k \too 0
$$
we get an exact sequence:
$$
L^1i^*i_*i^*E^k \too i^*E^{k+1} \too i^*E^k \too i^*i_*i^*E^k \too 0.
$$
Here,
the morphism $i^*E^k\to i^*i_*i^*E^k$ 
is an isomorphism by~\ref{iprop}~(9).
It follows that
 $i^*E^{k+1}$ is a quotient of the Artin sheaf 
$L^1i^*i_*i^*E^k\simeq E^k(-1),$
see \ref{iprop}~(10). Since a quotient of an Artin sheaf on
the commutative $\PP^1$ is obviously Artin again,
it follows that $i^*E^{k+1}$ is an Artin sheaf. 
The induction hypothesis implies that $E^{k+1}$ is an Artin sheaf. Hence
$E$, being an extension of Artin sheaves, is also Artin by~\ref{artin}~(5). 
\qed\medskip

\begin{definition}
A sheaf $E$ on $\ppgt\tau$ is called {\em $z$-torsion free} if either
of the following equivalent conditions hold {\rm(}equivalence is
proved in $\ref{iprop}~(2)${\rm)}

\newi1 the map $z:E\to E(1)$ is a monomorphism;

\newi2 $L^1i^*E=0$.
\end{definition}

\begin{lemma}\label{ztf}
\newi1 If $E$ is torsion free then $E$ is $z$-torsion free.

\newi2 For any coherent sheaf $E$ on $\ppgt\tau$ there exists
a unique subsheaf $F\subset E$ supported on $\plg$ such that $E/F$
is $z$-torsion free.
\end{lemma}
{\sl Proof:}
\newi1 Follows from \ref{iprop}~(5).

\newi2 First we will check the existence.
Consider the sequence:
 $F_k=\Ker(E\stackrel{z^k}\to E(k))\,,\,k=0,1,2,\ldots,$
of subseaves in $E$. Since the algebra $A^\tau$ is noetherian
the sequence $F_k$ stabilizes. Thus we have $F_{n+k}=F_n$
for some $n$ and for all $k\ge0$. The sheaf $F=F_n$ is clearly
supported on $\plg$. We claim that $E/F$ is $z$-torsion free.
To see this, note that the
equality: $F_n=F_{n+1}$ means that the composition
of the embedding: $E/F\into E(n)$ with the morphism:
$E(n)\stackrel {z\cdot}{\too} E(n+1)$ is an embedding. But this composition
can be factored through the morphism: $E/F\stackrel {z\cdot}{\too}
 (E/F)(1)$. Thus,
 $E/F$ is $z$-torsion free.

It remains to prove the uniqueness. Assume that $F'$ is a 
 subsheaf of $E$ supported
on $\plg$ and such that $E/F'$ is $z$-torsion free.
Then for any $k\in\Z$ we have an exact sequence
$$
\Ker \bigl(F' \stackrel{z^{n+k}}\too F'(n+k)\bigr) \into
\Ker \bigl(E \stackrel{z^{n+k}}\too E(n+k)\bigr) \to
\Ker \bigl(E/F' \stackrel{z^{n+k}}\too (E/F')(n+k)\bigr)
$$
It is clear that for $k\gg0$ the first term in the above
sequence coincides with $F'$, the second term coincides with
$F_{n+k}=F$ and the third term vanishes. Thus $F'=F$.
\qed\medskip
\nopagebreak

\subsection{Properties of $\coh(\PP^2_{\!_\G})$ for generic $\tau$}

The definition of "generic" parameters
$\tau$, due to [CBH], has been sketched in the Introduction.
The only property that will be used below is that,
for generic $\tau$, the algebra $\B^\tau$ has no
non-trivial finite-dimensional modules. 
In the context of deformed preprojective algebras,
an equivalent
propery 
has been proved in [CBH, \S7].

\begin{proposition}\label{mainiprop} Suppose that $\tau$ is
generic. Then:

\newi1 If $i^*E=0$ then $E=0$.

\newi2 If $\phi\in\Hom(E,F)$ and $i^*\phi$ is an epimorphism, then
$\phi$ is an epimorphism.

\newi3 If $\phi\in\Hom(E,F)$ and both $i^*\phi$ and $L^1i^*\phi$
are isomorphisms then so is $\phi$.

\newi4 If $\phi\in\Hom(E,F)$, $i^*\phi$ is a monomorphism and $L^1i^*F=0$
then $\phi$ is a monomorphism.

\newi5 A sheaf $E$ is locally free if and only if  $L^{>0}i^*E=0$ and $i^*E$ is
locally free.

\newi6 If $E\in\coh(\PP^2_{\!_\G})$ is torsion free, and the sheaf
$i^* E$ is locally free, then $E$ is locally free.
\end{proposition}
{\sl  Proof:}
\newi1 The  statement translated into the language of
$A^\tau$-modules reads: {\it If $M$ is a finitely generated graded
$A^\tau$-module such that $\dim_{_\C}(M/zM)<\infty$,
then $\dim_{_\C}(M)<\infty$.} To prove this, 
Note first that,
for any $i\geq 0$, the space $M_i$ is finite
dimensional, since $M$ is finitely generated.
Now, assume  $M/zM$ is finite dimensional. Hence, the multiplication map
$z:M_i\to M_{i+1}$ is surjective, for $i\gg0$.
It follows that, for $i\gg0$, the sequence:
$\dim M_i \geq \dim M_{i+1}\geq \ldots,$
stabilizes, hence the map $z:M_i\to M_{i+1}$
is an isomorphism. But for such an $i$, the endomorphisms
$z^{-1}x,z^{-1}y:M_i\to M_i$ provide $M_i$  with the structure of 
a finite
dimensional $\B^\tau$-module. For generic $\tau$, the
 algebra $\B^\tau$ has no finite-dimensional representations (because it
is Morita equivalent to the deformed preprojective algebra,
and the latter has no finite-dimensional representations by 
Theorem 7.7 of \cite{CBH}). Therefore $M_i=0$ for $i\gg0$, hence $M$
is finite dimensional.

\newi2 Assume that $\phi\in\Hom(E,F)$ is such that $i^*\phi$ is an
epimorphism. Since $i^*$ is right exact it follows that
$i^*\Coker\phi=0$, hence $\Coker\phi=0$ by $(1)$.

\newi3 It follows from $(2)$ that $\phi$ is an epimorphism.
Hence we have an exact sequence
$$
L^1i^*E   \quad\stackrel{L^1i^*\phi}\lra\quad  L^1i^*F  \quad\lra\quad
i^*\bigl(\Ker\phi\bigr)
 \quad\lra\quad i^*E \quad\stackrel{i^*\phi}\lra\quad i^*F
\quad\lra\quad 0.
$$
It follows that $i^*\Ker\phi=0$, hence $\Ker\phi=0$ by (1), that is
$\phi$ is a monomorphism.  Thus $\phi$ is an isomorphism.

\newi4 If  $L^1i^*F=0$ then  $i^*\Ker\phi=\Ker i^*\phi =  0$, hence
$\Ker\phi=0$ by $(1)$.

\newi5   If  $E$   is  locally   free  then   by   \ref{iprop}(4)  and
\ref{iprop}(6)  we have  $L^{>0}i^*E=0$  and $i^*E$  is locally  free.
Conversely,  assume  that  $L^{>0}i^*E=0$.  Then we  have  a  spectral
sequence
$$
E_2^{p,q}        =       L^{-q}i^*\CExt^p(E,\CO)       \Longrightarrow
\CExt^k(i^*E,\CO_{\plg}).
$$
It follows from \ref{iprop}~(3)  that the spectral sequence
degenerates in the second term.  Hence, if $i^*E$ is locally free,
we get
$$
E_2^{1,0}  =  i^*\CExt^1(E,\CO)  = 0\quad\text{and}\quad  E_2^{2,0}  =
i^*\CExt^2(E,\CO) = 0
$$
Thus, $E$ is locally free by  (1).

\newi6 Follows from (5) and \ref{iprop}~(5)
\qed\medskip

\section{Interpretation of quiver varieties}
\setcounter{equation}{0}

\subsection{From quiver data to a sheaf}\label{fromqtos}
Let $L^*$ be the dual of the tautological 2-dimensional $\G$-module.
We fix $V$ and $W$, finite dimensional $\Gamma$-modules, and a triple:
$$
(B,I,J) \in \Hom_{_\G}(V,V\otimes_{_\C} L^*) \;\bplus\;
\Hom_{_\G}(W,V)\; \bplus\; \Hom_{_\G}(V,W).
$$

Let 
$\{e_x,e_y\}$ be the
basis of $L^*$ dual to the basis $\{x,y\}$ of $L$. Then we can
consider $e_xx+e_yy \in L^*\gtimes L$  as an
element of
$$
 L^*\,\,\mbox{$\bigotimes$}_{_\CG}\; (L \oplus \triv)=
L^*\gtimes A_1 =
H^0\bigl(\PP^2_{\!_\G}\,,\,L^*\gtimes\CO(1)\bigr) =
\Hom\bigl(\CO(-1)\,,\,L^*\gtimes\CO\bigr),
$$
where $\triv$ stands for the trivial 1-dimensional $\G$-module.
Dually, viewing $e_xx+e_yy$ as an element of $L\gtimes L^*$,
we get a natural element in:
$$
(L\oplus \triv)\;\mbox{$\bigotimes$}_{_\CG}\; L^* =
A_1\gtimes L^* =
\Hom\bigl(L\gtimes\CO\,,\,\CO(1)\bigr)\,.
$$
Thus, we can define canonical  sheaf morphisms: $a = a_{_{B,I,J}}$
and $b=b_{_{B,I,J}}$,  by:
\begin{align}\label{theab}
&a = \left(\!\!
\begin{array}{c}B\cdot z - \id_V
\otimes(e_x\cdot x+e_y\cdot y)\\J\cdot z\end{array}
\!\!\right)\;:\;\;
V\gtimes\CO(-1) \too
\bigl((V\otimes_{_\C}L)\,\, \bplus\,\, W\bigr)\gtimes\CO\\
& 
b = \left(
B\cdot z - 
\id_V\otimes(e_x\cdot x+e_y\cdot y)\,,\,\,I\cdot z
\right)\;:\;\;
\bigl((V\otimes_{_\C}L)\,\, \bplus\,\, W\bigr)\gtimes\CO
\too V\gtimes\CO(1)
.\nonumber
\end{align}

From now on we use the canonical $\G$-module isomorphism:
$L^*\simeq L$, to identify $L^*$ with $L$. Thus, the maps 
(\ref{theab}) give 
the following morphisms of coherent
sheaves on $\PP^2_{\!_\G}$:
\begin{equation}\label{themonada}
V\gtimes\CO(-1)\;\; \stackrel{a_{_{B,I,J}}}\too \;\;
\bigl((V\otimes_{_\C}L)\,\, \bplus\,\,
 W\bigr)\gtimes\CO \;\;\stackrel{b_{_{B,I,J}}}\too\;\;
V\gtimes\CO(1)\,.
\end{equation}

\begin{definition}
A {\em monad} is a three-term complex, ${\mathcal{C}}$,
concentrated in degrees:
$-1,0,1,$ with the single non-zero cohomology group, $H^0({\mathcal{C}})$,
 referred to
as the {\em cohomology} of the monad.
\end{definition} 

\begin{proposition}\label{q-s}
If $a_{_{B,I,J}}$ and $b_{_{B,I,J}}$ are given by $(\ref{theab})$ 
with $(B,I,J)\in \bbm^\tau_{_\G}(V,W),$
cf. (\ref{bbm}),
then $(\ref{themonada})$ is a monad and its middle cohomology sheaf
$E$ admits a canonical framing $i^*E\cong W\gtimes\CO$ and,
moreover,
$H^1(\PP^2_{\!_\G},E(-1))\cong V$.
\end{proposition}

The proof will take the rest of this subsection. It will be 
divided into a sequence of lemmas. We begin with an obvious

\begin{lemma}\label{istarm}
The restriction of~$(\ref{themonada})$ to $\PP^1_{_\G}$ is
a complex which is
canonically quasi-isomorphic to $W\gtimes\CO$.
$\quad\square$
\end{lemma}

\begin{lemma}\label{ba0}
The triple $(B,I,J)$ satisfies the Moment Map Equation if
and only if ${b_{_{B,I,J}}\ccirc a_{_{B,I,J}}=0}.$ 
\end{lemma}
{\sl Proof:}
Straightforward computation shows that
$b_{_{B,I,J}}\ccirc a_{_{B,I,J}} = ([B,B] + IJ - \tau)\cdot z^2$.
\qed\medskip

From now on assume that $(B,I,J)$ satisfies the Moment Map
Equation. Then (\ref{themonada}) is a complex by
Lemma \ref{ba0}. We choose
the cohomological grading of this complex so that its middle term has
degree zero. Let $\CH^p$, $p=-1,0,1,$ denote the cohomology sheaves of
this complex.

\begin{lemma}\label{cohartin}
We have $i^*\CH^{-1}=i^*\CH^1=0$ and $i^*\CH^0$ is a subsheaf in
$W\gtimes\CO$. In particular, $\CH^{-1}$ and $\CH^1$ are Artin
sheaves. If $W=0$ then $\CH^0$ is Artin as well.
\end{lemma}
{\sl Proof:}
It follows from Lemma \ref{istarm} that we have a spectral sequence with the
second term 
$$
E^2_{p,q} = L^{-q}i^*\CH^p \Longrightarrow E^\infty_i = 
\begin{case}
W\gtimes\CO, & \text{if $i=0$}\\0, & \text{otherwise}\end{case}
$$
On the other hand~\ref{iprop}~(3)  implies that $E^2_{p,q}=0$, for
 $q\ne 0,1$, hence this spectral sequence
degenerates at the second term. Hence $i^*\CH^{-1}=i^*\CH^1=0$ and
$i^*\CH^0$ is a subsheaf in $W\gtimes\CO$. But
then~\ref{iprop}~(10) implies that $\CH^{-1}$ and $\CH^1$ are Artin
sheaves, and $\CH^0$ is also Artin whenever $W=0$.
\qed\medskip

\begin{lemma}\label{ain}
The map $a_{_{B,I,J}}$ is injective.
\end{lemma}
{\sl Proof:}
We 
note that the sheaf $\CH^{-1}=\Ker a$ is simultaneously an Artin sheaf
by \ref{cohartin} and locally free by \ref{lfprop}~(1). Hence it
vanishes by~\ref{artin}~(4). 
\qed\medskip

\begin{lemma}\label{bsur}
If $(B,I,J)$ is stable then $b_{_{B,I,J}}$ is surjective.
\end{lemma}
{\sl Proof:}
The sheaf $\CH^1=\Coker b$ is Artin by~\ref{cohartin}. Let
$\phi: V\gtimes
\CO(1)\onto\CH^1$ be the canonical projection and
let 
$$
V' = \Ker\Big[\phi(-1):\,H^0(\PP^2_{\!_\G}\,,\,
V\gtimes\CO) \to 
H^0(\PP^2_{\!_\G},\CH^1(-1))\Big].
$$
Note that $z$-multiplication gives an isomorphism:
 $\CH^1(-1)\iso\CH^1$
(see~\ref{cohartin} and \ref{iprop}~(10)). Hence the condition
$\phi\cdot b=0$ implies: $I(W)\subset V'$ and 
$B(V')\subset V'\otimes L$. Stability of $(B,I,J)$ then yields:
$V'=V$. Since $\phi$ is surjective and $\CH^1$ is Artin it follows
that $\CH^1=0$, that is, $b$ is surjective.
\qed\medskip

\begin{lemma}\label{h0tf}
The sheaf $\CH^0$ is torsion free.
\end{lemma}
{\sl Proof:}
Let $\CC$ denote the complex~(\ref{themonada}) and let $\CC^*$ denote the
dual complex (in the category $\qgr\lef(A^\tau)$
of sheaves of left modules). We have
a spectral sequence:
$$
E^2_{p,q} = \CExt^q\bigl(\CH^{-p}(\CC)\,,\,\CO\bigr) \Longrightarrow 
E^\infty_i = \CH^i(\CC^*).
$$
Since $\CH^{-1}=0$ by \ref{ain},
 it follows that the spectral sequence
degenerates at the third term. Moreover,
$\CExt^2(\CH^0(\CC),\CO)=0,$ and $\CExt^1(\CH^0(\CC),\CO)$ is a
quotient of the sheaf $\CH^1(\CC^*)$. On the other hand, the complex
$\CC^*$ is the complex, corresponding to the dual quiver data
$(B^*,J^*,I^*)$ (in the category of sheaves of left modules), hence by
\ref{cohartin} we have:
 $i^*\CH^1(\CC^*)=0$. Since $i^*$ is right exact
it follows that $i^*\CExt^1(\CH^0(\CC),\CO)=0$. Therefore,
$\CExt^1(\CH^0(\CC),\CO)$ is Artin by \ref{iprop}~(10). 

Now consider the spectral sequence
$$
E^2_{p,q} = \CExt^q\Bigl(\CExt^{-p}\bigl(\CH^0(\CC)\,,\,
\CO\bigr)\,,\,\CO\Bigr) \Longrightarrow 
E^\infty_i = 
\begin{case}\CH^0(\CC), & \text{if $i=0$}\\0, & \text{otherwise}\end{case}
$$
We have already
proved that $\CExt^2\bigl(\CH^0(\CC),\CO\bigr)=0,$ and that
$\CExt^1\bigl(\CH^0(\CC),\CO\bigr)$ is Artin. It follows that the spectral
sequence degenerates at the third term giving rise to a
short exact sequence
(see~\ref{artin}~(3)):
$$
0 \to \CH^0(\CC) \to \bigl(\CH^0(\CC)\bigr)^{**} \to
\CExt^2\Bigl(\CExt^1\bigl(\CH^0(\CC)\,,\,
\CO\bigr)\,,\,\CO\Bigr) \to 0.
$$
The middle sheaf here is locally free by \ref{lfprop}~(5), hence
$\CH^0(\CC)$ is torsion free by definition.~$\square$
\medskip

\noindent
{\bf Proof of proposition~\ref{q-s}.}\quad 
If $(B,I,J)\in \bbm^\tau_{_\G}(V,W)$ then~(\ref{themonada}) is a complex by
\ref{ba0}, which is left exact by \ref{ain} and right exact by
\ref{bsur}. Hence it is a monad. Moreover, its middle cohomology sheaf
$E$ is torsion free by \ref{h0tf} and admits a canonical framing by
\ref{istarm}. Finally, it is easy to see that
$H^1(\PP^2_{\!_\G},E(-1))\cong V$. \qed\medskip

Associating to any quiver data $(B,I,J)\in \bbm^\tau_{_\G}(V,W)$, the middle
cohomology sheaf of the corresponding monad~(\ref{themonada}) we
obtain a map:
 $\bbm^\tau_{_\G}(V,W) \to \bcm^\tau_{_\G}(V,W)$. It is clear that this
map is $G_{_\G}(V)$-equivariant. Indeed, any element
$g\in G_{_\G}(V)$ gives an isomorphism between the complex corresponding to  
a
quiver data $(B,I,J)$ and the
complex corresponding to the quiver data
$(gBg^{-1},gI,Jg^{-1})$. It follows that the corresponding middle
cohomology sheaves are isomorphic. Thus we obtain a well-defined map:
\begin{equation}\label{map_mm}
\FM_{_\G}^\tau(V,W) \too \bcm^\tau_{_\G}(V,W),
\end{equation}
We will show that  this map provides the bijection claimed in
Theorem~\ref{quiv-sheaf}.

\subsection{From a framed sheaf to quiver data}

We are going to study framed torsion free sheaves on $\PP^2_{\!_\G}$
using the Beilinson spectral sequence. We will need the following
lemma, cf. \cite{KKO}.

\begin{lemma}\label{cohe}
Let $E$ be a framed torsion free sheaf on $\PP^2_{\!_\G}$. We have

\newi1 $H^0(\PP^2_{\!_\G}\,,\, E(-1)) = H^0(\PP^2_{\!_\G}\,,\, E(-2)) = 0$;

\newi2 $H^2(\PP^2_{\!_\G}\,,\, E(-1)) = H^2(\PP^2_{\!_\G}\,,\, E(-2)) = 0$;

\newi3 $\Hom(T(-1)\,,\, E(-1)) = \Ext^2(T(-1)\,,\, E(-1)) = 0 $;

\newi4 $H^1(\PP^2_{\!_\G}\,,\, E(-1)) \simeq H^1 (\PP^2_{\!_\G}\,,\, E(-2))$.
\end{lemma}
{\sl Proof:}
\newi1 We have: $L^0i^*E= W\gtimes\CO$, and $L^1i^*E=0$,
by Proposition  \ref{iprop}~(5). Thus,
the exact sequence of Proposition  \ref{iprop}~(2)
reads
\begin{equation}\label{ie}
0 \too E(k-1) \too E(k) \too i_*L^0i^*E(k)= i_*(W\gtimes\CO(k)) \too 0.
\end{equation}
Since 
$H^0\bigl(\PP^2_{\!_\G}\,,\,i_*(W\gtimes\CO(k))\bigr) = 
H^0\bigl(\PP^1_{_\G}\,,\,W\gtimes\CO(k)\bigr) = 0,$
for all $k<0$, we get
$$
H^0(\PP^2_{\!_\G}\,,\,E(-1)) = H^0(\PP^2_{\!_\G}\,,\,E(-2)) = \dots =
H^0(\PP^2_{\!_\G}\,,\,E(-k))\quad,\quad \forall k>0\,.
$$
Since $E$ is torsion free it can be embedded into a
sheaf $\CO(n)^{\oplus m}$, by~\ref{lfprop}~(4)). Hence
$$
H^0(\PP^2_{\!_\G}\,,\,E(-k)) \subset 
H^0(\PP^2_{\!_\G}\,,\,\CO(n-k)^{\oplus m}) = 0\quad,\quad \forall k>n\,.
$$
It follows that:
$H^0(\PP^2_{\!_\G}\,,\,E(-1)) = H^0(\PP^2_{\!_\G}\,,\,E(-2)) = 0$.

\newi2 Similarly, since 
$H^1(\PP^2_{\!_\G}\,,\,i_*(W\gtimes\CO(k))) = 
H^1(\PP^1_{_\G}\,,\,W\gtimes\CO(k)) = 0,$
for all $k\ge -1,$ we get
$$
H^2(\PP^2_{\!_\G}\,,\,E(-2)) = H^2(\PP^2_{\!_\G}\,,\,E(-1)) = \dots =
H^2(\PP^2_{\!_\G}\,,\,E(k))\quad,\quad k> 0.
$$
But the ampleness implies that,
for $k$ sufficiently large, one has:
$H^2(\PP^2_{\!_\G}\,,\,E(k))=0$. Thus,
$H^2(\PP^2_{\!_\G}\,,\,E(-2)) = H^2(\PP^2_{\!_\G}\,,\,E(-1)) = 0$.

\newi3 The first of the sequences (\ref{deft}) and part (1) imply:
$\Hom(T(-1)\,,\,E(-1))=0$. Similarly, the second of the sequences
(\ref{deft}) and part  (2) imply: $\Ext^2(T(-1)\,,\,E(-1))=0$. 

\newi4 Follows from (\ref{ie}) for $k=-1$.
\qed\medskip

\begin{corollary} \label{monya}
Any framed torsion free sheaf $E$ can be represented
as the middle cohomology sheaf of a monad of the form: 
$$
0 \to \widetilde{V} \gtimes \CO(-1) \to \widetilde{V}' \gtimes \CO \to
\widetilde{V} \gtimes \CO(-1) \to 0;
$$
where $\widetilde{V} = H^1(\PP^2_{\!_\G}\,,\, E(-1))$ and $\widetilde{V}' = \Ext^1(T,E)$.
\end{corollary}
{\sl Proof.} Follows from Beilinson's spectral sequence,
see \S3.1, applied to $E(-1)$,
and the vanishing results of Lemma \ref{cohe}.
\qed\medskip

Now, fix  $E\in \bcm^\tau_{_\G}(V,W),$  and put 
$\widetilde{V} = H^1(\PP^2_{\!_\G}\,,\, E(-1))$ 
and $\widetilde{V}' = \Ext^1(T,E)$,
and fix a  monad as in Lemma \ref{monya}.
Choose a triple  $(B,I,J)\in \bbm^\tau_{_\G}(V,W)$,
and consider the corresponding monad (\ref{themonada})

\begin{lemma}\label{frmonad}
Any vector space isomorphism
$\varphi: \widetilde{V}\iso V$ 
can be uniquely extended to an isomorphism of monads:
$$
\begin{array}{lccccc}
\mbox{\footnotesize monad \ref{monya}:} &\widetilde{V} \gtimes \CO(-1) 
&\into  &\widetilde{V}' \gtimes \CO &\onto
&\widetilde{V} \gtimes \CO(-1)\\
&\downarrow\varphi &&\downarrow&&\downarrow\\
\mbox{\footnotesize monad \ref{themonada}:} & V\gtimes\CO(-1) &
\stackrel{a_{_{B,I,J}}}\into&
\bigl((V\otimes_{_\C}L) \bplus W\bigr)\gtimes\CO &\stackrel{b_{_{B,I,J}}}\onto&
V\gtimes\CO(1) \;,
\end{array}
$$
which
is compatible with framings. 
\end{lemma}
{\sl Proof:}
A repetition of the proof in \cite{KKO}, Theorem 6.7.
\qed\medskip

\begin{lemma}\label{bijstable}
If the monad in Corollary $\ref{monya}$ is isomorphic to a
complex~$(\ref{themonada})$ 
for some $(B,I,J)$, then the data $(B,I,J)$ satisfy both the Moment Map
Equation and the Stability Condition.
\end{lemma}
{\sl Proof:}
The first part follows immediately from~\ref{ba0}. Thus we have to
check the stability. Assume that the triple $(B,I,J)$ is not stable
and let $V'\subset V$ be a $\Gamma$-submodule such that $V'\ne V$,
$I(W)\subset V'$, and $B(V')\subset V'\otimes L$. Let $V''=V/V'$ and
let $\CC(V,W)$, $\CC(V',W)$ and $\CC(V'',0)$ be the complexes of the
form~(\ref{themonada}) given by the triples, induced by
$(B,I,J)$. Then we have
an exact sequence of complexes
$$
0 \to \CC(V',W) \to \CC(V,W) \to \CC(V'',0) \to 0,
$$
which gives rise to a long exact sequence of cohomology sheaves:
\begin{equation}\label{long}
\ldots\to \CH^0\CC(V'',0)\to \CH^1\CC(V',W) \to
\CH^1\CC(V,W)\to\CH^1\CC(V'',0) \to 0\,.
\end{equation}
Note that $\CH^{-1}\CC(V'',0)=0$ 
by~\ref{ain}. The cohomology $\CH^0\CC(V'',0)$
is both an Artin sheaf by~\ref{cohartin},
and a torsion free sheaf, by~\ref{h0tf}.
It follows that $\CH^0\CC(V'',0)=0$,
by \ref{tfprop} and
\ref{artin}~(4).
Assume: $\CH^1\CC(V'',0)=0$. Then the complex
$\CC(V'',0)(-1)$ is quasi-isomorphic to zero, hence
the hypercohomology spectral sequence with $E_2^{p,q}$-term:
$\dis E_2^{p,q}=H^p\bigl(\PP^2_{\!_\G}\,,\,
\CH^q\CC(V'',0)(-1)\bigr)$
would converge to zero. This sequence, however, clearly
converges to the vector space $V''$.  Therefore,
 $\CH^1\CC(V'',0)\neq 0$.
But
then the long exact sequence (\ref{long})
would force: $\CH^1\CC(V,W)\neq 0,$ which
is a contradiction, because $\CC(V,W)$ is a monad.
The contradiction implies that $V''=0$.
\qed\medskip

\begin{proposition}\label{mapiso}
The map: $\FM_{_\G}^\tau(V,W)\to \bcm^\tau_{_\G}(V,W)$ defined in 
(\ref{map_mm}) is a bijection.
\end{proposition}
{\sl Proof:}
Surjectivity of the map follows from~\ref{monya}, \ref{frmonad} and
\ref{bijstable}. To check injectivity, note that it follows
from~\ref{frmonad} that the quiver data $(B,I,J)$ are uniquely
determined by a torsion free framed sheaf $E$ up to 
isomorphism $H^1(\PP^2_{\!_\G},E(-1))\cong V$, that is up to the
action of
the group $G_{_\G}(V)$.
\qed\medskip

\noindent
Proposition~\ref{mapiso} completes the proof of
Theorem~\ref{quiv-sheaf}. 

We finish this section with a few remarks.

\begin{definition}
A triple $(B,I,J)$ is {\em costable} if for any $\Gamma$-submodule
$V'\subset V$ such that $J(V')=0$ and $B(V')\subset V'\otimes L$ we
have: $V'=0$.
\end{definition}

Note that a triple $(B,I,J)$ is costable if and only if  the dual triple
$(B^*,J^*,I^*)$ is stable.
\begin{proposition}[\cite{VV}]
The framed sheaf corresponding to a (stable)
quiver data\linebreak $(B,I,J)\in
\bbm^\tau_{_\G}(V,W)$ is locally
free if and only if  the triple $(B,I,J)$ is {\em costable}.
\end{proposition}
{\sl Proof:}
Repeating the arguments in the proof
of~\ref{h0tf} we see that for the middle cohomology sheaf $\CH^0$ 
of the complex~(\ref{themonada}) we have
$\CExt^2(\CH^0,\CO)=0$, while $\CExt^1(\CH^0,\CO)=0$ if and only if  the complex
corresponding to the dual triple is exact at the right term. But
\ref{bsur} and \ref{bijstable} imply that this holds if and only if  the dual
triple is stable.
\qed\medskip

\begin{remark}
In general, locally-free framed sheaves form an open dense  subset in 
$\bcm^\tau_{_\G}(V,W)$  since the Costability Condition is open. 
However, for a generic value of the parameter $\tau$ every stable
quiver data is  automatically costable, i.e. all torsion free 
framed sheaves are automatically locally free, cf.~\ref{mainiprop}~(6).
$\quad\lozenge$ 
\end{remark}

\section{Proof of the Crawley-Boevey and Holland Conjecture}
\setcounter{equation}{0}

Throughout this section we assume $\tau$ to be generic, though some
results remain true for an arbitrary $\tau$. 

Recall that  we have defined in the Introduction the algebra 
$\B^\tau=A^\tau/(z-1)\cdot A^\tau,$ where
$z$ is the degree one central variable in the algebra $A^\tau$.
Explicitly, we have: $\,\B^\tau=\C\langle x,y\rangle \# \Gamma\big/
\langle\!\langle [y,x] - \tau \rangle\!\rangle.\,$
The standard grading on the algebra $A^\tau$ induces a canonical
increasing filtration: $\CG=\B^\tau_0\subset \B^\tau_1\subset\B^\tau_2
\subset\ldots,$ on $\B^\tau$
 such that $\gr(\B^\tau)$,
the associated graded algebra, has finite
homological dimension.
Thus, it follows from a general result due to Quillen [Q], that
the assignment: $R \mapsto R\gtimes\B^\tau$
induces an isomorphism: $K(\Gamma)\iso K(\B^\tau)$
of the corresponding Grothendieck groups of projective modules.
 Let $[N]$ denote the class of a
$\B^\tau$-module $N$ in $K(\B^\tau)\cong K(\Gamma)$.

Recall that by Proposition \ref{easy},
to any 1-dimensional class $R\in K(\G)$,
one can canonically attach (isomorphism classes of)
$\G$-modules $W$ and $V$, such that in $K(\G)$ we have:
$R= [W]+ [V]\cdot([L]-2[\triv])$ and, moreover,
 $\dim W=1$ and  $V$ does not
contain the regular representation as a submodule.
The goal of this section is to prove the following

\begin{theorem}\label{lbl}
Given  a class $R\in K(\G)$ such that $\dim R=1$, let
$V,W$ be $\G$-modules attached to $R$ in Proposition \ref{easy}.
Then, there exists
a natural bijection:
$$
\left\lbrace
\begin{array}{l}\mbox{Isomorphism classes of finitely generated 
projective}\\
\mbox{$\B^\tau$-modules $N$ such that 
$[N]=R$ in $K(\B^\tau)=K(\Gamma)$}
\end{array}
\right\rbrace
\;\;\simeq\;\;
\bigsqcup_{k=0}^\infty\; \bcm_{_\G}^\tau(V\oplus \CG^{\oplus k}\,,\,W)
\,.
$$
\end{theorem}

This theorem together with Theorem~\ref{quiv-sheaf} 
yields Theorem~\ref{main_intro}. 

\subsection{From sheaves on $\PP^2_{\!_\G}$ to projective 
$\B^\tau$-modules}\label{fromstom}

Let $\modd(\B^\tau)$ denote  the category
of finitely generated (right) $\B^\tau$-modules.
There is a natural "{\it open restriction}" functor
$\,j^*: \gr(A^\tau)\too\modd(\B^\tau)\,,\, M \mapsto
 M/(z-1)\cdot M.$ 

It will be convenient  to use an equivalent definition
of the algebra $\B^\tau$ that will make 
the  open restriction functor $j^*$ manifestly {\it exact}.
 Namely, let $A^\tau[z^{-1}]$ denote 
 the localization of the
algebra $A^\tau=\bplus_k\;A^\tau_k$ with respect to $z$, and
$A^\tau[z^{-1}]_0,$ the degree zero component of the localized
algebra. We have:
\[
\B^\tau \;\;\simeq\;\; A^\tau[z^{-1}]_0\; \;=\;\; \dlim A^\tau_k\,,
\]
where the direct 
limit is taken with respect to the embeddings: 
$A^\tau_k\stackrel{z\cdot}{\too} A^\tau_{k+1},$ induced by
multiplication by $z$. Using this formula one can rewrite
the functor  $j^*$ in the form:
$$
j^*:\; M=\bplus_k\,M_k \quad\mapsto\quad j^* M = \dlim M_k\,,
$$
where the limit is taken with respect to the embeddings:
$M_k\stackrel{z\cdot}{\too}
 M_{k+1},$ induced
by the $z$-action.

\begin{lemma}\label{jstar}
\newi1 The functor $j^*$ factors through the category
$\qgr(A^\tau)=\coh(\ppgt\tau)$.

\newi2 The functor $j^*$ is exact.

\newi3 The functor $j^*:\coh(\ppgt\tau)\to\modd(\B^\tau)$ commutes
with the dualization, i.e. for any coherent sheaf $E$ on $\ppgt\tau$
we have $j^*(E^*) = \Hom_{\modd(\B^\tau)}(j^*E,\B^\tau)$.

\newi4 For any $\Gamma$-module $R$ we have: 
$j^*(R\gtimes\CO(k))=R\gtimes\B^\tau$.
In particular, we have $[j^*(R\gtimes\CO(k))] = [R]\in K(\Gamma)$. 

\newi5 If $j^*E=0$ then the sheaf $E$ is supported on $\PP^1_{_\G}$.
In particular, $r(E)=0$.
\end{lemma}
{\sl Proof:}
\newi1 It suffices to check that, for any finite
dimensional graded $A^\tau$-module $M$, we have:
 $j^*M=0$. But this is clear, because
in this case, for $k\gg0$, one has: $M_k=0$, hence $\dlim M_k=0$.

\newi2 It follows from the exactness of the direct limit.

\newi3 Let $M$ be the graded $A^\tau$-module corresponding to a
sheaf~$E$. Then we have
$$
\Hom_{\modd(\B^\tau)}(j^*E\,,\,\B^\tau) = 
\Hom_{\modd(\B^\tau)}(\dlim_k M_k\,,\, \B^\tau).
$$
The dual sheaf $E^*$ corresponds to the
$A^\tau$-module $\,\bplus_{\!k}\,
\Hom_{\coh(\ppgt\tau)}\bigl(E,\CO(k)\bigr)$,
by definition. Hence, we obtain
$$
\begin{array}{l}
\displaystyle
j^*(E^*) = \dlim_k\, \Hom_{\coh(\ppgt\tau)}(E\,,\,\CO(k)) \\
\displaystyle\qquad\qquad\qquad
= \dlim_k\, \Hom_{\qgr(A^\tau)}
\left(\mathop{\bplus}\limits_{n=0}^\infty M_n\,,\,A^\tau(k)\right)
= \dlim_{m,k}\, \Hom_{\gr(A^\tau)}
\left(\mathop{\bplus}\limits_{n=m}^\infty M_n\,,\,A^\tau(k)\right)\\
\qquad\qquad\qquad
=\dlim\limits_m\,\Hom_{\gr(A^\tau)}
\left(\mathop{\bplus}\limits_{n=m}^\infty M_n\,,\,
\dlim\limits_k\,A^\tau(k)\right)
=
\Hom_{\modd(\B^\tau)}(\dlim\limits_n M_n\,,\,\B^\tau).
\end{array}
$$
This is precisely what we need.

\newi4 We have:
$\,\dis
j^*(M\gtimes\CO(k)) = \dlim_l \left(M\gtimes A^\tau_{k+l}\right) =
M\gtimes\bigl(\dlim_l A^\tau_{k+l}\bigr) = M\gtimes \B^\tau.
$

\newi5 Let $M$ be the graded $A^\tau$-module corresponding to a
sheaf~$E$. Then $j^*E=0$ implies that for any $m\in M$,
there exists $n\gg 0$ such that $mz^n=0$.
 Since $M$ is finitely generated we conclude
that $Mz^n=0$ for some $n\gg 0$. Let $M_k = \Ker z^k\subset M$
and put $E_k=\pi(M_k)$. This gives a filtration on $E$. Finally,
for each $k\in \Z$, the element
$z$ annihilates the quotient $M_k/M_{k-1}$, hence
$E_k/E_{k-1}=i_*F_k$ for some coherent sheaves $F_k$ on $\plg$.
\qed\medskip

\begin{proposition}\label{s-m}
If $E\in \bcm^\tau_{_\G}(V,W)$ is a framed torsion free sheaf then
$N:=j^*E$ is a projective $\B^\tau$-module with 
\begin{equation}\label{nvw}
[N]=[W] + [V\otimes L] - 2[V]\in K(\B^\tau) = K(\Gamma).
\end{equation} 
\end{proposition}
{\sl Proof:}
If $E$ is a torsion free sheaf, one can find
an embedding: $E\into \CO(n)^{\oplus m}$,
 for some $n$ and $m$. Applying $j^*$ we obtain an
embedding: $N\into (\B^\tau)^{\oplus m}$. On the other hand by
\cite{CBH} Theorem~0.4,
for generic  $\tau$, the global homological
dimension of the algebra $\B^\tau$
 equals~1. It follows that any submodule of a free
$\B^\tau$-module is projective. Thus $N$ is projective. 

It remains to compute the class of $N$ in $K(\B^\tau)$. To this end
we use the monadic description of torsion free sheaves
provided by Corollary \ref{monya} and Lemma \ref{frmonad}.
Writing $E$ as the cohomology of the monad
corresponding to a triple $(B,I,J)\in \bbm^\tau_{_\G}(V,W)$
and using \ref{jstar}~(1), (4), we find:
$$
\begin{array}{l}
[N] = [j^*E] = 
[j^*\bigl((V\otimes L \,\,\bplus\,\, W)\gtimes\CO\bigr)] 
\smallskip\\\qquad\qquad
- [j^*\bigl(V\gtimes \CO(-1)\bigr)] 
- [j^*\bigl(V\gtimes \CO(1)\bigr)] =
[W] + [V\otimes L] - 2[V]. \qquad\square
\end{array}
$$
\medskip

Note, that if $R=\CG^{\oplus k}$ is a multiple of the regular representation of
$\Gamma$ then we have an isomorphism of $\Gamma$-modules: 
$R\otimes L\cong R\oplus R$, hence,
$[R\otimes L] - 2[R]=0$. Therefore, given
a 1-dimensional
$\Gamma$-module $W$ and a $\Gamma$-module $V$ that {\it does not contain
the regular representation as a submodule}, 
we see from Proposition
\ref{s-m} that the assignment: $E \mapsto j^*E$ gives 
a map
\begin{equation}\label{eqq}
\begin{array}{rcl}
\bigsqcup\limits_{k=0}^\infty\; \bcm^\tau_{_\G}(V\oplus\CG^{\oplus k}\,,
\,W)
& \stackrel{j^*}{\too} &
\Big\{
{\begin{array}{c}
{\text{projective $\B^\tau$-modules $N$ such that}}\\
{[N]=[W]+[V]\cdot([L]-2[\triv])}
\end{array}}
\Big\}\,.
\end{array}
\end{equation}
(here $[N]\in K(\B^\tau)$ is treated as a class in
$K(\G)$ via the isomorphism:
$ K(\B^\tau)\simeq K(\G)$, as above).
We will prove below that (\ref{eqq}) is a bijection. This
will imply Theorem \ref{lbl}.

\subsection{Extending $\B^\tau$-modules to sheaves on $\PP^2_{\!_\G}$}
In this subsection
we  show that, given  a projective
$\B^\tau$-module
$N$, there exists an essentially unique (up to isomorphism) way
to extend $N$ to a framed torsion free 
 sheaf $E$ on $\PP^2_{\!_\G}$ such that $j^*E\cong N$.

Recall first that 
the standard grading on the algebra $A^\tau$ induces a canonical
increasing filtration: $\CG=\B^\tau_0\subset \B^\tau_1\subset\B^\tau_2
\subset\ldots,$ on the algebra $ \B^\tau=A^\tau/(z-1)A^\tau$.
Given a $ \B^\tau$-module $N$, we say that
an increasing filtration $\,\{N_k\}\,$ on  $N$ is 
{\em compatible} with the canonical filtration on $\B^\tau$ if, for
all $k$, $l,$ we have $N_k\cdot\B^\tau_l\subset N_{k+l}$. 
The filtration $\,\{N_k\}\,$
 is said to be {\em finitely generated} if  $\,\dim
N_k<\infty\,,\,\forall i,$ 
 and there exists $k$ such that $N_k\cdot\B^\tau_l =
N_{k+l}$ for all $l\ge 0$. The filtration 
is called {\em exhausting}
if $N =\bigcup_k N_k$. Finally, two
filtrations: $\,\{N_k\}\,$  and $\,\{N'_k\}\,$  on $N$
are called {\em equivalent}, for all
 $k\gg0$, we have:  $N_k = N'_k$.

\begin{proposition}\label{extj1}
The set of $z$-torsion free coherent sheaves $E$ on $\ppgt\tau$ such
that $j^*E\cong N$ is in  bijection with the set of equivalence  
classes of finitely generated increasing exhausting filtrations $\,\{N_k\}\,$
on $N,$ compatible with the canonical filtration of the
algebra~$\B^\tau$.  
\end{proposition}
{\sl Proof:}
\newi1 If $M=\bplus\,M_k$ is the graded $A^\tau$-module corresponding to a
$z$-torsion free sheaf $E$ then, for $k\gg0$, the $z$-multiplication
map: 
$M_k\stackrel{z\cdot}{\too} M_{k+1}$ is injective.
Hence the images of $\,\{M_k\}_{k\ge 0}\,$ form an
increasing filtration on $\,\dlim M_k = j^*E.$
This filtration is clearly
finitely generated, exhausting and compatible with the canonical
filtration of $\B^\tau$. 

Conversely, assume that $N$ is a $\B^\tau$-module with a finitely
generated increasing exhausting filtration $\,\{N_k\}\,$ compatible with the
canonical filtration of $\B^\tau$. Then the graded vector space 
$\bplus_k\, N_k$ admits the structure of a graded
$A^\tau$-module by means of the standard Rees construction.
Specifically, we let $x,y \in A^\tau$ act
on $N=\bplus_k\, N_k$ via the maps
 $\,x,y:N_k\to N_{k+1}$ induced by the same named elements
of $\B^\tau$, and we let  $z$ act as the tautological embedding:
$N_k\into N_{k+1}$. It is clear that the graded $A^\tau$-module 
thus defined is
finitely generated and the corresponding coherent sheaf $E$ on
$\PP^2_{\!_\G}$ is $z$-torsion free. Moreover, it is easy to show that
equivalent filtrations give rise to isomorphic coherent
sheaves. Finally, it is clear  that the construction
of this paragraph is 
inverse to that of the preceding one.
\qed\medskip

\begin{proposition}\label{extj2}
If $E$ and $E'$ are $z$-torsion free sheaves on $\ppgt\tau,$
and $\phi:j^*E\to j^*E'$ is a morphism of $\B^\tau$-modules then
there exists $n\ge0$ and a
morphism $\tilde\phi:E\to E'(n),$ such that
$j^*\tilde\phi=\phi$.
\end{proposition}
{\sl Proof:}
Let $N=j^*E$, $N'=j^*E'$ and let $\,\{N_k\}\,$,
$\,\{N'_k\}\,$ be the corresponding
finitely generated filtrations on $N$ and $N'$. Since the filtration
$\,\{N_k\}\,$ is finitely generated it follows that there exists $n\ge 0$ such
that $\phi(N_k)\subset N'_{k+n},$ for all $k\gg0$. Hence $\phi$ gives a
morphism of graded $A^\tau$-modules: $\bplus_{\!k}\; N_k \too
\bplus_{\!k}\; N'_{k+n}$
or, equivalently, a morphism of coherent sheaves $\tilde\phi:E\to
E'(n)$. It is clear that $j^*\tilde\phi=\phi$.
\qed\medskip

\begin{lemma}\label{extj3}
For any coherent sheaf $E$ on $\PP^2_{\!_\G}$ we have:
$\;
r(E) = \dim_{_\C}[j^*E],
\;$ 
where $\dim_{_\C}:K(\Gamma)\to\Z$ is the linear function given by
$\dim_{_\C}([R]) = \dim_{_\C} R.$
\end{lemma}
{\sl Proof:}
Note that both RHS and LHS are linear functions on
$K(\coh(\PP^2_{\!_\G}))$, see~\ref{rri} and \ref{jstar}~(2). Thus it
suffices to verify the equality only for $E=R\gtimes\CO(i)$,
see~\ref{k0}. This has been  done in \ref{rri} and \ref{jstar}~(4).
\qed\medskip

\subsection{Proof of  bijectivity}

From now untill the end of 
this section we fix a class
$R\in K(\G)$ such that $\dim R=1$. By Proposition \ref{easy}, in $K(\G)$ we can write:
$R= [W]+ [V]\cdot([L]-2[\triv]),$
for certain uniquely determined (isomorphism classes of)
$\G$-modules $W$ and $V$, such that $\dim W=1$ and 
such that $V$ does not
contain the regular representation as a submodule.
With  $W$ and $V$ as above, we have
\begin{proposition}\label{surj}
If $N$ is a projective finitely generated $\B^\tau$-module such that
$[N]=R$,
then there exists a framed locally free sheaf 
$\CE\in \bcm^\tau_{_\G}(V\oplus\CG^{\oplus k}\,,\,W)$ 
such that $j^*\CE\cong N$.
\end{proposition}
{\sl Proof:}
Choose an arbitrary finitely generated increasing exhausting
filtration on $N$ compatible with canonical filtration of $\B^\tau$
and let $E$ be the corresponding $z$-torsion free sheaf on
$\PP^2_{\!_\G}$ such that $j^*E\cong N$ (see \ref{extj1}). Now,
for 
the sheaf $E^{**}$, by~\ref{jstar}~(3) we have:
$$
j^*(E^{**}) = \Hom_{\B^\tau}\bigl(j^*(E^*)\,,\,\B^\tau\bigr) =
\Hom_{\B^\tau}\bigl(\Hom_{\B^\tau}(N,\B^\tau)\,,\,\B^\tau\bigr) = N,
$$
since $N$ is projective. On the other hand,
 $E^{**}$ is a locally free sheaf
by~\ref{lfprop}~(5). Hence $i^*E^{**}$ is locally free
by~\ref{iprop}~(6). Moreover, by~\ref{iprop}~(7) and \ref{extj3},
we get
$$
r(i^*E^{**}) = r(E^{**}) = \dim_{_\C}[j^*E^{**}] =
\dim_{_\C}[N] = \dim_{_\C}[R] = 1
$$
Hence by~\ref{plg}~(3) we have:
$i^*E^{**}\cong W\gtimes\CO(n),$ for a 1-dimensional
$\Gamma$-module $W$ and some $n\in\Z$.

Let $\CE=E^{**}(-n)$. Then $\CE$ is a locally free framed sheaf on
$\PP^2_{\!_\G}$, hence $\CE\in \bcm_{_\G}(V',W)^\tau$,
 for a certain
$\Gamma$-module $V'$. On the other hand it is clear that
$$
j^*\CE = j^*(E^{**}(-n)) \cong j^*E^{**} \cong N.
$$
Hence by \ref{s-m} we have:
$\,
[R] = [W] + [V\otimes L] -2[V],
\,$
hence Lemma \ref{easy} yields:
$V'\cong V\oplus\CG^{\oplus k}$, moreover,
$V$ and $W$ are $\G$-modules corresponding to the
class $[R]$ in the sense of Lemma \ref{easy}.
\qed\medskip

\begin{proposition}\label{inj}
Let $\CE\in \bcm^\tau_{_\G}(V\oplus\CG^{\oplus k}\,,\,W)$ and
$\CE'\in \bcm^\tau_{_\G}(V\oplus\CG^{\oplus k'}\,,\,W)$ be locally free
sheaves such that $j^*\CE\cong j^*\CE'$. Then,
$k=k'$, and  $\CE\cong\CE'$.
\end{proposition}
{\sl Proof:}
An isomorphism $j^*\CE\cong j^*\CE'$ gives by~\ref{extj2} a morphism
$\phi:\CE\to\CE'(n)$ such that $j^*\phi$ is an isomorphism. Let
$K=\Ker\phi$, $C=\Coker\phi$. Then by~\ref{jstar}~(2) we have
$j^*K=j^*C=0$, hence by~\ref{jstar}~(5) both $K$ and $C$ are supported
on $\plg$. On the other hand, a locally free
sheaf contains no sheaves supported on $\plg$,  by~\ref{ztf}~(2) and
(1).
It follows that  $K=0$. Thus we
have a short exact sequence:
\begin{equation}\label{defc}
0 \too \CE \stackrel\phi\too \CE'(n) \too C \too 0.
\end{equation}
Applying the functor $i^*$ we get a short exact sequence:
$$
0 \too L^1i^*C \too W\gtimes\CO \stackrel{i^*\phi}\too 
W\gtimes\CO(n) \too i^*C \too 0.
$$
Since $\dim_{_\C} W=1$ it follows that either $i^*\phi=0$ or $i^*C$ is an
Artin sheaf. 

If $i^*C$ is Artin then by~\ref{supi}~(3) the sheaf $C$ is Artin
as well.
On the other hand, applying $\CHom(-,\CO)$ to~(\ref{defc})
we get an exact sequence
$$
\CExt^1(\CE,\CO) \too \CExt^2(C,\CO) \too \CExt^2(\CE'(n),\CO)
$$
But $\CE$ and $\CE'$ are locally free, hence $\CExt^2(C,\CO)=0$,
hence $C=0$ by~\ref{artin}~(4). Thus $\phi$ is an isomorphism.
It follows that $i^*\phi$ is an isomorphism, hence 
$W\gtimes\CO\cong W\gtimes\CO(n)$, hence $n=0$ and
$\CE\cong\CE'$.

If $i^*\phi=0$ it follows that $\phi$ factors through
the embedding: $\CE'(n-1) \stackrel {z\cdot}\into \CE'(n)$.
Repeating this argument for $n$ being replaced by $n-1,n-2,\ldots,$
we obtain that either $\CE\cong\CE'$ or
there exists an embedding: $\CE\into\CE'(n)$, for arbitrarily
small $n\in\Z$. The latter  is impossible by \ref{resolution} and
\ref{lfprop}~(4). 
\qed\medskip

Now, Proposition~\ref{surj} gives surjectivity
of the map (\ref{eqq}), and Proposition~\ref{inj} gives  injectivity
of (\ref{eqq}). Hence, this map is bijective, and Theorem
\ref{lbl} follows. 

\section{Appendix A: Graded preprojective algebra.}
\setcounter{equation}{0}

In this section we define  a graded version $\bp$ of the deformed
preprojective algebra introduced in  \cite{CBH}.  Let $Q$ be a quiver,
i.e. an oriented graph  with 
vertex set $\Vx$. For
any (oriented) edge $a\in Q$, we write $in(a)  = j$, $out(a) = i,
$ if $a: i \to
j$. Let $\overline{Q}$ be the  double of $Q$,
obtained   by  adding a
reverse edge $a^*: j \to i$ for every edge $a: i \to j$ in $Q$.

Let  $\bp_0=\bplus_{\!v\in\Vx}\;\C$  be the
direct sum of $|\Vx|$ copies of the field $\C$, a
commutative  semisimple  $\C$-algebra. 
For $v\in \Vx$, 
we write $e_v\in \bp_{\!_0}$ for the projector on the $v$-th copy
(an idempotent). We define a $\bp_{\!_0}$-bimodule $\bp_{\!_1}$ by
the formula
$$
\bp_{\!_1} = \Big(\mathop{\bplus}\nolimits_{a \in \overline{Q}}
\; \C \cdot a \Big)
\;\bigoplus\;
\bp_{\!_0}\;.
$$
Here, in the first summand, for
 $a\in\overline{Q}$, with
$in(a)  = j$, $out(a) = i,
$ we put: $e_i  a = a e_j =
a$,  and all  other products:  $e_k a$,  $a e_s$ 
 are set equal to zero.  
Let $f$ denote the canonical generator of the second
summand in $\bp_{\!_1}$ corresponding to the element $1\in\bp_{\!_0}$.
We put $f_i=e_i\cdot f= f\cdot e_i,$ so that $f=\sum f_i$.

The $\bp_{\!_0}$-bimodule $\bp_{\!_1}$ gives rise to the tensor algebra
$T^\bullet_{_{\bp_{\!_0}}}(\bp_{\!_1}) = \bplus_{n \geq 0}\;
T^n_{_{\bp_0}}\bp_{\!_1}$, where $T^n_{_{\bp_0}}\bp_{\!_1}
 = \bp_{\!_1} \otimes_{_{\bp_{\!_0}}} \ldots
\otimes_{_{\bp_{\!_0}}}  \bp_{\!_1}$ is the $n$-fold tensor product. Note that,
since the product is taken over $\bp_{\!_0}$, for any  two arrows $a, a'
\in \overline{Q}$, in $T^\bullet_{_{\bp_{\!_0}}}(\bp_{\!_1})$
 we have: $a \cdot a' = 0$ unless $in(a') = out(a)$.

\bigskip 
\noindent
\textbf{Definition.}   
Choose an  element  $\tau \in  \bp_{\!_0}$, $\tau  =
\sum_{i=0}^n \tau_i  e_i$.  The \textit{graded  deformed preprojective
algebra}, $\bp^\tau =\bp^\tau(Q)$, is defined as
 $\bp^\tau  =T^\bullet_{_{\bp_{\!_0}}}(\bp_{\!_1})/  \langle\!\langle
 R\rangle\!\rangle,$  a quotient of
 the  tensor  algebra $T^\bullet_{_{\bp_{\!_0}}}(\bp_{\!_1})$  by the  
two-sided  ideal   generated  by the
$\bp_{\!_0}$-bimodule  $R  
\subset  \bp_{\!_1}  \otimes_{_{\bp_{\!_0}}}  \bp_{\!_1}$  formed
by the
following quadratic relations:

(a) $f_i\cdot
 a = a\cdot f_j$$\;,$\quad if $a: i \to j$ is an arrow in $\overline{Q}$

(b)  $  \displaystyle{\sum_{\{a \in  Q  \;|\; out(a)  =  i\}}}  a\cdot
a^*  
\;\,-
\displaystyle{\sum_{\{ a \in Q \;|\; in(a) = i\}}} a^*\cdot a = \tau_i\cdot f_i^2,
\qquad \forall i \in \Vx$.

\noindent
\textbf{Koszul complex.}  By the  definition the algebra $\bp^\tau $ is
quadratic (see Appendix  B). Therefore one can write  its right Koszul
complex $\K^\bullet\bp^\tau$
(see  Apendix B for  the definition, and also  \cite{BGS}). In
our particular case it boils down to
$$
0 \to \K^3 \otimes_{_{\bp_{\!_0}}} \bp^\tau (-3) \too  \K^2 \otimes_{_{\bp_{\!_0}}}
\bp^\tau (-2) \too \K^1 \otimes_{_{\bp_{\!_0}}} \bp^\tau (-1) \too \bp^\tau 
\too \bp_{\!_0} \to 0
$$
where $\K^i=\K^i\bp^\tau $ are $\bp_{\!_0}$-bimodules given by:

(a) $\K^1\bp^\tau  = \bp_{\!_1}$;

(b) $\K^2\bp^\tau =R \subset \bp_{\!_1} \otimes_{_{\bp_{\!_0}}} 
\bp_{\!_1},$\quad is the submodule of
generating relations;

(c) $\K^3\bp^\tau  =\bigl(\K^2 \otimes_{_{\bp_{\!_0}}} \bp_{\!_1}\bigr)\,
\bigcap\, \bigl(\bp_{\!_1}
\otimes_{_{\bp_{\!_0}}} \K^2\bigr)\,
\subset\, T^3_{_{\bp_{\!_0}}}\bp_{\!_1}$ is
a $\bp_{\!_0}$-bimodule that can be shown to have a single generator:
$$
\tau \cdot f \otimes f \otimes f\, +\, \sum_{a \in Q} a \otimes f \otimes
a^* - f \otimes a \otimes a^* - a^* \otimes f \otimes a + f \otimes
a^* \otimes a - a \otimes a^* \otimes f + a^* \otimes a \otimes f
$$
The differentials of the Koszul complex are given by  restricting
 the map:
$$
T_{_{\bp_{\!_0}}}^n\bp_{\!_1} \otimes_{_{\bp_{\!_0}}} \bp^\tau  (-1) \to
T_{_{\bp_{\!_0}}}^{n-1}\bp_{\!_1}
 \otimes_{_{\bp_{\!_0}}} \bp^\tau ; \qquad
(v_1 \otimes \ldots \otimes v_n) \otimes x \mapsto
(v_1 \otimes \ldots \otimes v_{n-1}) \otimes v_n x.
$$

\bigskip 

The algebra $\,^!  \bp^\tau $, see Appendix B for the general
definition of the dual quadratic algebra,
is generated by the $\bp_{\!_0}$-bimodule
$\,^! \bp_{\!_1}$ which is spanned over $\C$ by two collections of elements:
$\,\{b\}_{b\in\overline{Q}},\,$ and
$\,\{r_i\}_{i\in\Vx},\,$ subject to the following relations 
$$
\begin{array}{ll}
\mbox{(a) $b\cdot r_i + r_j\cdot  b =
0$} 
&\mbox{if $b: i \to j$ is an arrow in $\overline{Q}$;}\\
\mbox{(b) $b_1\cdot  b_2 = 0$ }
&\mbox{unless $b_1 \in Q\;\&\; b_2 =
b_1^*,$ or $b_2 \in Q\;\&\; b_1 = b_2^*$;}\\
\mbox{(c) $\tau_i\cdot  b\cdot  b^* = r_i^2$} &\mbox{if $b
\in Q$ and $in(b) = i$}\\
\mbox{(d) $\tau_i\cdot  b^*\cdot  b = r_i^2$} &\mbox{if $b
\in Q$ and $out(b) = i$}\\
\mbox{(e) $b^*_1 \cdot b_1 = b_2\cdot 
b^*_2$}  &\mbox{if $b_1, b_2 \in Q$ and $in(b_2) = i =
out(b_1).$}
\end{array}
$$
(relation (e)
does not follow from (c) and (d) if and only if  $\tau_i = 0$).
One can check that the relations above imply ${}^! \bp^{\tau}_i = 0$,
for all $i \geq 4$.

It is known, see [GMT],  that,  for  any  quiver   $Q$,  the
corresponding  algebra  $\bp^\tau =\bp^\tau (Q)$  is \textit{Koszul}.
However,  $\bp^\tau $ is  noetherian and has polynomial
growth if and only if the underlying  graph  $Q$ is
either of  affine  or of finite Dynkin ADE-type.

We assume that $Q$ is an affine Dynkin graph,
hence it is associated, by means of McKay correspondence,
to a finite subgroup $\G\subset SL_2(\C)$.
We write $Q=Q(\G)$, and
let $R_i$ be the simple $\G$-module corresponding to a vertex
$i\in\Vx$.
Given $\tau\in \ZZ(\CG)$, let $\tau_i\in\C$ be the complex number
such that $\tau$ acts in $R_i$ as the scalar operator:
$\tau_i\cdot\id_{_{R_i}}$. This way we identify $\tau$
with the element $\sum_i\,\tau_i\cdot e_i\in\bp_{\!_0}$,
still to be denoted by  $\tau$.
With this understood, one proves as in [CBH]:
\begin{proposition}
Let $Q(\G)$ be the graph of affine ADE-type arising
from a finite group $\G\subset SL_2(\C)$ via the
MacKay corresondence. Then

\vi The algebras $A^\tau$ and $\bp^\tau =\bp^\tau(Q)$ are Morita
equivalent. In particular,

\vii The category $\gr(A^\tau)$
is equivalent to  the category $\gr(\bp^\tau )$,  resp. category
$\qgr(A^\tau)$ is equivalent to $\qgr(\bp^\tau )$.$\quad\square$
\end{proposition}

\begin{definition}
The deformed preprojective algebra of the quiver $Q$ is defined
as the quotient algebra
$$
\Pi^\tau\, :=\, \bp^\tau(Q)/\langle\!\langle f-\tau\rangle\!\rangle
\,=\,
\bp^\tau(Q)/\langle\!\langle f_i-\tau_i\rangle\!\rangle_{i\in\Vx}
$$
\end{definition}

Let $v\in \Vx$ be  a vertex of the graph  $Q=Q(\G)$,
and let $Q\fin$ be the graph obtained from $Q$ by removing `$v$'.
The vertex `$v$' is said to be a {\it special} vertex of
$Q$ if $Q\fin$ is a Dynkin graph of finite type and
the graph $Q$ is the  extended affine graph for $Q\fin$
(so that `$v$' becomes the extended vertex).
From Theorem \ref{main_intro} we deduce the following
generalization of the {\bf{Crawley-Boevey \& Holland conjecture}},
cf. [BLB, Example 5.7]:
\begin{corollary}\label{main}
Let $Q(\G)$ be the McKay graph of $\G$,
and $v\in \Vx$ a special vertex of $Q(\G)$. Let $\tau$ be generic,
and
$R$ be the 1-dimensional (simple) representation of $\G$ 
corresponding to the vertex $v$. Then, there exists
a natural bijection:
$$\left\lbrace
\begin{array}{l}\mbox{Isomorphism classes of finitely generated 
projective}\\
\mbox{$\Pi^\tau$-modules $N$ such that 
$[N]=[R]$ in $K(\Pi^\tau)=K(\Gamma)$}
\end{array}
\right\rbrace
\;\;\simeq\;\;
\bigsqcup_{k=0}^\infty\; \FM_{_\G}^\tau(V\oplus\CG^{\oplus k}\,,\,W)\,.\,\square
$$
\end{corollary}
Here we have used the  natural isomorphism:
$K(\Pi^\tau)\simeq K(\Gamma)$, see [Q].\nopagebreak
\section{Appendix B: Algebraic generalities.}
\setcounter{equation}{0}
\subsection{Linear algebra  over a semisimple algebra \cite[\S2.7]{BGS}}
Let $A_0$ be a finite dimensional semisimple $\C$-algebra.

Recall that for any \textit{left} $\anol$-module $V$ the
space $V^* = Hom_{\anol\mbox{\tiny{\sf{-mod}}}}(V, \anol)$ 
can be given the structure of a
\textit{right} $\anol$-module via the assignment: $(fa) (v) =
(f(v))a$. Similarly, for any \textit{right} $\anol$-module $W$ the
space ${}^*W = Hom_{\mbox{\tiny{\sf{mod-}}}\anol}(W, 
\anol)$ can be given the structure of
a \textit{left} $\anol$-module via the assignment: $(ag)(w) = a
(g(w))$. For finitely generated left $A_0$-modules
$V, W,$ the canonical  evaluation maps: $V \to
{}^*(V^*)$ and $W \to (^*W)^*$  are isomorphisms.

For  an $\anol$-bimodule  $V$, both  $V^*$ and  ${}^*V$  are bimodules
defined as follows: $(af)(v) =  f(va)$ and $(ga)(v) = g(av)$, $\forall
g \in {}^*V$, $f \in V^*$, $v \in V$, $a \in \anol$.

All standard results of linear  algebra over a field can be generalized
in an appropriate way to $\anol$-modules (e.g. $W^* \otimes V^* \simeq
(V \otimes W)^*$) if we take  all tensor products over $\anol$. We will
use  these generalizations  freely, referring  the interested  reader to
\cite[\S2.7]{BGS}.

\medskip
Fix  an algebra $A= \bigoplus_{n \geq 0} A_n$, and put $X=\proj A$.
The following result has been proved in [AZ, Theorem 8.1(3)]:

\begin{proposition}\label{sek}If $A$ is strongly regular
of dimension $d$, then one has:
\begin{equation*}
H^p(X,\CO(i))=\begin{case} A_i, & \text{if $p=0$ and $i\ge0$}\\
{}^*  \!  A_{-i-d-1},  &   \text{if  $p=d$  and  $i\le-d-1$}\\  0,  &
\text{otherwise}
\end{case}
\end{equation*}
where ${}^* \! A_{-i-d-1}$ is the  dual of the $A_0$-bimodule
$A_{-i-d-1}$ in the sense explained above.
\end{proposition}
{\sl Sketch of Proof.}
Let $M$ be a graded module and $F$ the corresponding sheaf.
For any fixed $i\geq 0$,  consider the right $A$-module
$\;\bplus_{n = -\infty}^{\infty}\; H^i (F(n)).$
It was shown in [AZ, Proposition 7.2(2)] that, for $i \geq 1$, one has
$$
\mathop{\bplus}_{n = -\infty}^{\infty}\; H^i (F(n)) \;\;\simeq\;\;
\lim_{m \to \infty}\;\Bigl( \mathop{\bplus}_{n = - \infty}^{\infty}\; 
\Ext^i\bigl(A/A_{\geq m}\,,\,
M(n)\bigr)\Bigr)\quad,\quad
A_{\geq m}=\bplus_{k\geq m}\;A_k\,.
$$
Now take $M = A$ and apply the above formula. Since $A/A_{\geq m}$
is finite dimensional, the Gorenstein condition implies that the LHS
above vanishes for $i= 1,\ldots, d-1$. 
This yields the cohomology vanishing part of the formula of the Proposition.
Other claims follow easily using Serre duality.\qed

\subsection{Koszul and co-Koszul algebras}
From now on we assume that $A = \bigoplus_{n \geq 0} A_n$ is a
positively graded $\C$-algebra such that all graded components
$A_n$ are finite
dimensional over $\C$.

\medskip
\noindent
\textbf{Definition.}
An algebra $A = \bigoplus_{n \geq 0} A_n$ is called
\textit{quadratic} if

\noindent
$\bullet\quad$ $\anol$ is a finite-dimensional semisimple $\C$-algebra;

\noindent
$\bullet\quad$ the $\anol$-bimodule $A_1$ generates $A$ over $\anol$;

\noindent
$\bullet\quad$ the relations ideal is generated by the subspace of
  quadratic relations $R \subset A_1 \otimes_{\anol} A_1$.
\medskip

Given a quadratic  algebra $A= \bigoplus_{n \geq 0} A_n$,
we can  represent $A$  as $T_{\anol}A_1/\langle\!\langle
R\rangle\!\rangle$,
the  quotient of  the tensor
algebra 
by the ideal $\langle\!\langle R\rangle\!\rangle$ generated  by the space
of quadratic relations $R \subset A_1 \otimes_{\anol} A_1$.

Define $A^!$, the \textit{left dual} of $A$, to be the quadratic  algebra:
$T_{\anol} (A_1^*)/\langle\!\langle R^{\perp}\rangle\!\rangle,$
 with $R^{\perp} \subset 
A_1^*
\otimes_{\anol} A_1^* = (A_1 \otimes_{\anol} A_1)^*$ being the
annihilator of $R$. Analogously, the \textit{right dual} ${}^!A$ is
defined as the algebra: $T_{\anol}({}^*\!A_1)/\langle\!\langle {}^{\perp}R
\rangle\!\rangle,$ with ${}^{\perp}R
\subset {}^*\! A_1 \otimes_{\anol} {}^*\! A_1 = {}^*(A_1 \otimes_{\anol}
A_1)$.

The \textit{right Koszul complex}  $\K^\bullet\!A$ is 
a complex of the form,
see e.g., [BGS],[Ma]:
$$
\ldots \stackrel{d}\longrightarrow ({}^!A_3)^* \otimes_{\anol} A(-3)
\stackrel{d}\longrightarrow ({}^!A_2)^* \otimes_{\anol} A(-2)
\stackrel{d}\longrightarrow ({}^!A_1)^* \otimes_{\anol} A(-1)
\to A \to \anol \to 0
$$
where the differential
`$d$' is  defined  as  follows.  Observe  that:  $\,({}^!   A_i)^*
\otimes_{\anol} A = $\linebreak
$ Hom_{\anol\mbox{\tiny{\sf{-mod}}}}({}^!A_i, A)$. 
Under the canonical
isomorphism $Hom_{\mbox{\tiny{\sf{mod-}}}\anol
}(A_1, A_1)  = A_1 \otimes_{\anol} {}^*\!A_1$,
let $\,\id_{A_1} = \sum v_{\alpha} \otimes \check{v}_{\alpha}$. Then for $f
\in  Hom_{\anol\mbox{\tiny{\sf{-mod}}}}({}^!A_{i+1}, A)$
  and  $a \in  {}^!  A_i$ we  set
$df(a) = \sum v_{\alpha}\cdot f(\check{v}_{\alpha}\cdot  a)$.
One can check  that this formula indeed defines  a complex. 
One can  also define a left  Koszul complex. It  is known, cf.  for example
\cite{BGS},  that  the  exactness  of  the right  Koszul  complex  is
equivalent to the exactness of the left Koszul complex.

Similarly, there is a natural differential on the
space $\K_\bullet A={}^!A\otimes_{\anol}  A,$ see [Ma],
making it into a complex, called the (right) {\it co-Koszul complex} of $A$.

\begin{definition}\label{cokos}\vi
 A quadratic ring $A$ is called \textit{Koszul} if
its (right) Koszul complex, $\K^\bullet\!A$
 has the only non-trivial cohomology in degree
zero.

\vii  A quadratic ring $A$ is called \textit{co-Koszul} of degree
$d$,  if
its (right) co-Koszul complex, $\K_\bullet{A}$
 has the only non-trivial cohomology in degree $d$ and, moreover,
$H^d(\K_\bullet{A})\simeq A_0(d)$.
\end{definition}

The conditions in the Proposition below are the basic conditions
that allow us to start out with the non-commutative
geometry as discussed in \S2.

\begin{proposition}\label{kosz}
Let $A= \bigoplus_{n \geq 0} A_n$ be a quadratic algebra
with ${}^!\!A=\bplus_{\!k}\;{}^!\!A_k$. Assume
 that:

$\bullet\quad$ ${}^!\!A$
has no non-zero graded components in degrees $k>d$,
and ${}^!\!A_d\simeq A_0$;

$\bullet\quad$The  algebra $A$ is Noetherian, and has
polynomial
growth;

$\bullet\quad$The  algebra $A$ is both
Koszul and co-Koszul (of degree $d$).

Then $A$ is strongly regular of dimension $d$
 in the sense of Definition \ref{reg}.
 \end{proposition}
{\sl Proof:}
To prove that the global dimension, $gl.dim(A)$, is finite we apply
\cite{Hu} and conclude that  $gl.  dim (A)$ equals
the  minimal length of projective resolution  for $A_0$.
The Koszul complex, if exact, provides such a minimal resolution.
Thus, for a Koszul algebra $A$, the global dimension equals the 
number of non-zero graded components of the algebra ${}^!\!A$.
Since $\dim_{_\C}({}^!\!A)$ is finite,
we conclude that: $gl.dim(A)<\infty$.

Notice next that an obvious canonical isomorphism:
 ${}^!\!A\otimes_{\anol}A \simeq
\Hom_{A}\bigl(({}^!\!A)^*\otimes_{\anol}A
\,,\,A\bigr),\,$ gives an isomorphism of complexes:
$\K_\bullet{A}\simeq\Hom_{A}\bigl(\K^\bullet\!A\,,\,A\bigr).$
It follows that, for a Koszul algebra $A$, the
complex $\K_\bullet{A}$ computes the $\Ext$-groups:
$\Ext^\bullet_A(A_0, A)$. Thus, $A$ is co-Koszul of degree $d$
if and only if it is
Gorenstein with parameters $(d,d)$.\quad \qed\medskip

\begin{remark}
  One shows similarly that if  $A$  is Koszul  of global  dimension $d$  and
  Gorenstein, then it is  Gorenstein with parameters $(d,d)$, coKoszul
  of degree $d$  and, moreover, the dual algebra ${}^!\!A$  is Frobenius of
  index $d$.$\quad\lozenge$
\end{remark}\bigskip

Fix   $A= \bigoplus_{n \geq 0} A_n$, an algebra satisfying the
conditions of Proposition \ref{kosz}, and put $X=\proj{A}$.
A key role in our study of sheaves on  $X$ is played by 
\medskip

\noindent
{\!\bf{
Beilinson spectral sequence (\cite{KKO}):}}
For any sheaf $E$ on $X$ there is a spectral sequence with the first term
$$
E_1^{p,q} = \Ext^q(Q_{-p}(p),E)\otimes_{\anol}\CO(-p)
\Longrightarrow E_\infty^i =
\begin{case} E, & \text{for $i=0$}\\ 0, & \text{otherwise}\end{case}
$$
where $p=-d,\dots,0$, and $Q_{-p}$ is the sheaf on $X$ corresponding
to the cohomology $\widetilde Q_{-p}$ of the truncated Koszul complex:
$$
0 \to A \to A_1^! \otimes_{\anol} A(1) \to \dots \to
A_{-p}^! \otimes_{\anol} A(-p) \to \widetilde Q_{-p} \to 0
$$
(here $A(n)$ stands for the algebra $A$ with the grading being
shifted by
$n$).
Equivalently, $\widetilde Q_{-p}$ can be described as follows
$$
0 \to \widetilde Q_{-p} \to A_{1-p}^! \otimes_{\anol} A(1-p) \to \dots \to
A_{d+1}^! \otimes_{\anol} A(d+1) \to \anol(d+1) \to 0.
$$
Note that all the sheaves $Q_p$ are naturally $A_0$-bimodules,
hence the tensor product in the $E_1$-term makes sense.
\subsection{Special  case: ${\mathbf{A=A}^{\boldsymbol{\tau}}}$}

Recall that:
$A^\tau =  \Bigl(\bigl(TL[z]\bigr)\# \G\Bigr)
\Big/  
\big\langle
\!\big\langle u\cd v - v\cd u - \om(u,v)\cdot\tau  z^2 \big\rangle
\!\big\rangle_{u,v\in L},\;$
see Definition \ref{Atau}. Write $\Lambda L^*$ for the
exterior algebra of the 2-dimensional vector space $L^*$,
and let $\Lambda L^*{}^{\,}\dot{\otimes}{}^{\,}\C[\xi]$ 
be the {\it super-}tensor product of $\Lambda L^*$
with the polynomial algebra in an {\it odd} variable $\xi$
of degree 1.
Thus, by definition, for any $v\in L^*\subset \Lambda L^*,$
in $\Lambda L^*{}^{\,}\dot{\otimes}{}^{\,}\C[\xi]$
we have: $v\cdot\xi=-\xi\cdot v$.
We will view $\Lambda L^*$ as a subalgebra in $
\Lambda L^*{}^{\,}\dot{\otimes}{}^{\,}\C[\xi]$.
Further, let $\om\in \Lambda^2L^*$ be the
element corresponding to the symplectic form on $L$.

\begin{proposition}
\label{kos-prep}
The algebra $A^\tau$ is a Noetherian algebra of polynomial
growth. Moreover, $A^\tau$ is both Koszul and co-Koszul (of degree 3),
and we have:
$${}^! \! A^\tau\;=\;
\Bigl(\bigl(\Lambda L^*{}^{\,}\dot{\otimes}{}^{\,}\C[\xi]\bigr)\# \G\Bigr)
\Big/  
\big\langle
\!\big\langle\om-\xi^2\cdot\tau\big\rangle
\!\big\rangle\;.
$$
\end{proposition}
{\sl Proof:}  
For $\tau = 0$  we have:
$A^\tau= \C[x,y,z]\#\G$. In this case
all the claims are easy and  follow e.g.,
 from  \cite{GMT}. Next, one checks
that the relations defining $A^\tau$ are 3-self-concordant in the 
sense of \cite{Dr}. Hence the graded components of $A^\tau$ are 
isomorphic to those of $A^0$ as vector spaces (cf. \cite{Dr}).
This implies that $A^\tau$ has polynomial
growth, for any $\tau$. Furthermore, it follows from the
Drinfeld's result that we
may view the family of Koszul complexes  $\K^\bullet{A^\tau}$,
resp. $\K_\bullet{A^\tau}$, as a family
of varying (with $\tau$)
 differentials on the  Koszul 
complex for $A^0$. Since the differential
for $\tau = 0$ has a single
non-trivial cohomology, the same is true for 
all values of $\tau$
close enough to zero. However, the algebras $A^\tau$ and $A^{\alpha\cdot
 \tau}$
are isomorphic for any $\alpha \in \C^*$.
Thus, $A^\tau$ is both Koszul and co-Koszul, for any $\tau$.

The expression for ${}^!\!A^\tau$ is obtained by a
direct calculation. It shows, in particular, that the algebra
${}^!\!A^\tau$ has non-vanishing graded components in degrees $i=0,1,2,3$ only.
Thus, $d=gl. dim(A^\tau)=3$. In particular,
formula (\ref{cohoi}) follows from Proposition \ref{sek}.
\qed\medskip

Let $R_i$, $i \in \Vx,$ be a
complete collection of the isomorphism classes of 
simple modules over
$\anol = \CG$, and  $\CR_i =
\pi(R_i \otimes_{_\CG} A^\tau)$ the coherent sheaf corresponding to
the graded right $A^\tau$-module $R_i \otimes_{_\CG} A^\tau$. Since
$R_i$ is a direct summand in $\anol$, the sheaf $\CR_i$ is a
direct summand of $\CO$. In particular, for any
$i\in\Vx$, the sheaf $\CR_i$ is locally
free in the sense of Definition \ref{free}.

\begin{remark}
Consider collection $\,\{
\CR_0,\dots,\CR_n,\dots,\CR_0(d-1),\dots,\CR_n(d-1)\}\,$
of sheaves on $\PP^2_{\!_\G}
=\proj A^\tau$. It follows from the explicit form of the cohomology
of the sheaves $\CO(n)$  that the above collection is
a {\em strong exceptional collection}, see [Ru].
 The Beilinson spectral sequence
of a sheaf $E$ can be considered as a decomposition of $E$ with respect
to the above exceptional collection.$\quad\lozenge$
\end{remark}


\begin{proposition}\label{resolution_App}
Any coherent sheaf $E$ on $\ppg$  admits a resolution of the form
$$
0 \too \mathop{\bplus}\nolimits_{i\in\Vx}\,\, 
V^d_i \otimes_{\C} \CR_i (k-d) \too \dots
  \too \mathop{\bplus}\nolimits_{i\in\Vx}\,\,
 V^0_i \otimes_{\C} \CR_i (k) \too E \too 0
$$
where $V^d_i,\dots,V^0_i$ are certain complex vector spaces.
\end{proposition}
{\sl Proof:} Consider the Beilinson spectral sequence of the sheaf
$E(n)$. By ampleness, for $n\gg0$ all higher $\Ext$-groups in the spectral
sequence vanish and only the $q=0$ row of it will be non-trivial. This
gives a resolution of  type (\ref{oresol}) for the sheaf $E(n)$.
Now tensor it with $\CO(-n)$. \qed\medskip

\section{Appendix C: Minuscule classes}
\setcounter{equation}{0}

\nc{\TC}{{\widehat{C}}}
\nc{\TP}{{{\sf{\hat{P}}}}}
\nc{\TD}{{\widehat{Delta}}}
\nc{\ttheta}{{\hat\theta}}
\nc{\tctheta}{{{\hat\theta}^\vee}}
\nc{\ctheta}{{\theta^\vee}}

\newtheorem{reform}{Reformulation}

The goal of this section is to prove Proposition \ref{easy}.
Thus, we fix
$\G$, a finite subgroup in $SL_2$ and let $Q(\G)$
denote the corresponding affine Dynkin graph.
We identify the set $I$ of vertices of  $Q(\G)$
with simple roots of the affine root system
associated to $Q(\G)$, and  write $\omega_i$, resp. $R_i$,
for the fundamental weight, resp.  irreducible $\G$-module,
corresponding to the vertex $i\in I$.
Thus, the $\omega_i$'s form a basis of the
weight lattice $\TP$ (of the affine root system), and the $R_i$'s
form a basis of $K(\G)$. The correspondence: $
\omega_i\longleftrightarrow R_i$
yields an isomorphism of lattices: $\TP\cong K(\G)$.

Let $L$ denote the tautological two dimensional representation
of $\Gamma$ and let $\triv$ be the trivial (one dimensional)
representation of $\Gamma$.
The map $K(\G)\to K(\G)\,,\,
[V]\mapsto$\linebreak
$[V]\otimes([L]-2[\triv])$ gets
identified under the isomorphism $\TP\cong K(\G)$
with the Cartan operator $\TC:\TP\to\TP$,
i.e., the linear map given (in the basis $\{\omega_i\}$) by
 the Cartan matrix.

Let $\TP^*=\Hom(\TP,\Z)$ be the coroot lattice,
and $\TC^*: \TP^*\to \TP^*$ the dual Cartan operator.
Note that $\dim\Ker\TC=1$ and $\dim\Ker\TC^*=1$, since
our root system is of affine type. Let
$\ttheta\in\TP$ and $\tctheta\in\TP^*$ denote the minimal
positive elements in $\Ker\TC$ and in $\Ker\TC^*,$
respectively. 
In the other words, $\ttheta$ and $\tctheta$ are
the minimal positive imaginary root and coroot, respectively. 
The class in $K(\G)$ of the regular representation of $\G$ gets identified
with $\ttheta\in \TP$, while the dimension function $\dim:
K(\G)\to\Z$
gets identified with the element
$\tctheta\in \TP^*$ considered as a function $\TP\to\Z$.

We see that the  isomorphism classes of one dimensional
$\Gamma$-modules are in bijection with
{\it special}  vertices of $Q(\G)$, i.e.,
the vertices  $i\in I$
such that $\dim[R_i]=1$, or equivalently $\tctheta(\omega_i)=1$.
A vertex $i\in I$ is known to be special if and only if
the graph $Q(\G)$ is obtained from a finite Dynkin
graph $Q\fin$ with vertex set $I\smallsetminus\{i\}$ by adding the vertex $i$.

Now Proposition \ref{easy} can be reformulated as follows.

\begin{lemma}\label{reformulation1}
For any $\omega\in\TP$ such that $\tctheta(\omega)=1$ there exists
a uniquely determined pair $(i,\omega_0)$, where $i\in I$ is a special
vertex and $\omega_0\in \TP$ is a dominant weight, such that the weight
$\omega_0-\theta$ is not dominant and
$$
\omega = \omega_i + \TC(\omega_0).
$$
\end{lemma}

Since $\Ker\TC=\Z\ttheta$ it follows that for any $\omega\in\Im\TC$
there exists a unique $\omega_0\in\TP$ such that $\omega_0$ is
dominant, $\omega_0-\theta$ is not dominant and $\omega=\TC(\omega_0)$.
Hence Lemma \ref{reformulation1} is equivalent to

\begin{lemma}\label{reformulation2}
For any $\omega\in\TP$ such that $\tctheta(\omega)=1$ there exists
a uniquely determined special vertex $i$ such that
$$
\omega - \omega_i \in \Im\TC.
$$
\end{lemma}

 From now on we fix some special vertex $v\in I$ and let
$I\fin=I\smallsetminus\{v\}$ be the vetex set of the
corresponding Dynkin graph $Q\fin$
of finite type. Let $\BP$ be the weight lattice of $Q\fin$, let
$C:\BP\to \BP$ be its Cartan operator, and let $\theta\in \BP$
and $\ctheta\in \BP^*$ denote the maximal root and coroot
in the root and coroot system of $Q\fin$ respectively.

The decomposition $I = I\fin\sqcup\{v\}$ gives rise
to the direct sum decompositions $\TP = \BP\oplus\Z\cdot\omega_{v}$
and $\TP^* = \BP^*\oplus\Z\cdot \alpha_v^\vee$, where $\alpha_v^\vee$ is the
simple coroot of $Q(\G)$ corresponding to the vertex $v$.
It is well known that we have
$$
\tctheta = (\ctheta,\alpha_v^\vee),\qquad
\ttheta  = (\theta,\omega_{v}),\quad\mbox{and}\quad
\TC = \left(
\begin{array}{ccc}C & , & \bullet \\ \bullet & , & \bullet \end{array}
\right).
$$
It follows that the projection ${\mathtt{pr}}:\TP= \BP\oplus\Z\cdot\omega_{v}
\to \BP$ gives rise
to the following
isomorphisms
\begin{equation}\label{xiiso}
\begin{array}{rcl}
\Big\{\omega\in\TP\ |\ \tctheta(\omega)=1 \Big\} & \cong & \BP \\
\Big\{i\in I\ |\ \tctheta(\omega_i)=1 \Big\}    & \cong &
\{v\} \bigsqcup \Big\{i\in I\fin\ |\ \ctheta(\omega_i)=1 \Big\}
\end{array}
\end{equation} 

It follows from (\ref{xiiso}) that if $\tctheta(\omega)=1$ then the condition
\begin{equation}\label{cond2}
\omega = \omega_i + \TC(\omega'),\quad\mbox{\sl where $i$ is a special 
vertex in $I$}
\end{equation}
is equivalent to the condition
$\,
{\mathtt{pr}}(\omega) = {\mathtt{pr}}(\omega_i) + {\mathtt{pr}}(\TC(\omega'))
.$
On the other hand, it is clear that
$\omega' = \omega'' + \alpha_v^\vee(\omega_{v})\cdot\omega_{v}$,
for some $\omega''\in \BP$. Hence
$\,
{\mathtt{pr}}(\TC(\omega')) = {\mathtt{pr}}(\TC(\omega'')) = C(\omega'').
$ 
Thus  condition (\ref{cond2}) is equivalent to
$$
{\mathtt{pr}}(\omega) = {\mathtt{pr}}(\omega_i) + C(\omega''),\quad
\mbox{\sl where $i$ is a special vertex in $I\fin$, or $\omega_i=0$}\,.
$$
Thus we obtain the following reformulation of Lemma 
\ref{reformulation2}.

\begin{lemma}\label{reformulation3}
For any $\omega\in \BP$ there exists a uniquely determined pair
$(i,\omega'')$, where $i$ is either $v$ or a special vertex
of $I\fin$ and $\omega''\in \BP$ is a weight, such that
$$
\omega = \omega_i + C(\omega'').
$$
\end{lemma}

{\sl Proof of Lemma \ref{reformulation3}:}
The image of the Cartan operator $C$ is the root sublattice
$\QQ\subset \BP$. On the other hand, it is well known (see \cite{Bou},
\S2, Ex.~5) that each coset in $\BP/\QQ$ contains a unique minuscule
weight or zero. Finally, according to \cite{Bou}, \S1, Ex.~24
the set of minuscule weights coincides with the set of fundamental
weights of special vertices of $I\fin$.
$\qquad\Box$

\footnotesize{

}

\footnotesize{
{\bf V.B.}: Caltech, Mathematics 253-37,
Caltech,
Pasadena CA 91125, USA;\\
\hphantom{x}\quad\, {\tt baranovs@caltech.edu}}

\footnotesize{
{\bf V.G.}: Department of Mathematics, University of Chicago,
Chicago IL
60637, USA;\\
\hphantom{x}\quad\, {\tt ginzburg@math.uchicago.edu}}

\footnotesize{
{\bf A.K.}: Institute for Information Transmission Problems, Russia;\\
\hphantom{x}\quad\, {\tt sasha@kuznetsov.mccme.ru}}

\end{document}